\documentclass[aps,prd,twocolumn,showpacs,showkeys,preprintnumbers,amsmath,amssymb,floatfix,longbibliography,nofootinbib]{revtex4-2}
\pdfoutput=1
 
\usepackage{amsmath,latexsym}
\usepackage{xcolor}
\usepackage[colorlinks=true, urlcolor=blue, linkcolor=blue, citecolor=blue]{hyperref}
\usepackage{etoolbox}
\usepackage{breqn}

\usepackage{graphicx}
\usepackage{dcolumn}
\usepackage{bm}

\def\beq{\begin{equation}}
\def\eeq{\end{equation}}
\def\beqn{\begin{eqnarray}}
\def\eeqn{\end{eqnarray}}
\def\eeqno#1{\label{#1}\end{equation}}

\newcommand{\be}{\begin{equation}}
\newcommand{\ee}{\end{equation}}
\newcommand{\bea}{\begin{eqnarray}}
\newcommand{\eea}{\end{eqnarray}}

\makeglossary
\makeindex

\begin{document}

\preprint{JGM/015-The Isomorphism of H4 and E8}
\title[The Isomorphism of $H_4$ and $E_8$]{The Isomorphism of $H_4$ and $E_8$}
\url{https://www.TheoryOfEverything.org}

\author{J Gregory Moxness}
\homepage{https://www.TheoryOfEverything.org/theToE}
\email[mailto:jgmoxness@TheoryOfEverything.org]{}
\affiliation{TheoryOfEverything.org}

\date{October 30,2023}

\begin{abstract}
This paper gives an explicit isomorphic mapping from the 240 real $\mathbb{R}^{8}$ roots of the $E_8$ Gosset $4_{21}$ 8-polytope to two golden ratio scaled copies of the 120 root $H_4$ 600-cell quaternion 4-polytope using a traceless 8$\times$8 rotation matrix $\mathbb{U}$ with palindromic characteristic polynomial coefficients and a unitary form $e^{\text {i$\mathbb{U}$}}$. It also shows the inverse map from a single $H_4$ 600-cell to $E_8$ using a 4D$\hookrightarrow$8D chiral left$\leftrightarrow$right mapping function, $ \varphi$ scaling, and $\mathbb{U}^{-1}$. This approach shows that there are actually four copies of each 600-cell living within $E_8$ in the form of chiral $H_{4L}$$\oplus$$\varphi H_{4L}$$\oplus$$H_{4R}$$\oplus$$\varphi H_{4R}$ roots. In addition, it demonstrates a quaternion Weyl orbit construction of $H_4$-based 4-polytopes that provides an explicit mapping between $E_8$ and four copies of the tri-rectified Coxeter-Dynkin diagram of $H_4$, namely the 120-cell of order 600. Taking advantage of this property promises to open the door to as yet unexplored $E_8$-based Grand Unified Theories or GUTs.
\end{abstract}


\maketitle

\section{Introduction}

\begin{figure}[!ht]
\center
\includegraphics[width=200pt]{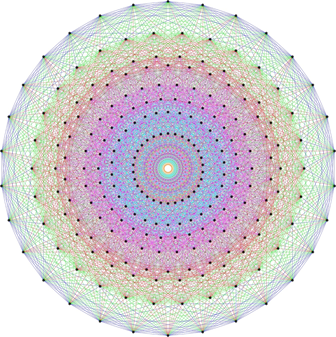}
\caption{\label{fig:E8Petrie} $E_8$ $4_{21}$ Petrie projection}
\end{figure}

Fig. \ref{fig:E8Petrie} is the Petrie projection of the Gosset $4_{21}$ 8-polytope derived from the Split Real Even (SRE) form of the $E_8$ Lie group with unimodular lattice in $\mathbb{R}^{8}$. It has 240 vertices and 6,720 edges of 8-dimensional (8D) length $\sqrt{2}$. $E_8$ is the largest of the exceptional simple Lie algebras, groups, lattices, and polytopes related to octonions ($\mathbb{O}$), (8,4) Hamming codes, and 3-qubit (8 basis state) Hadamard matrix gates. An important and related higher dimensional structure is the $\mathbb{R}^{24}$ ($\mathbb{C}^{12}$) Leech lattice ($\Lambda_{24}\supset$$E_8$$\oplus$$E_8$$\oplus$$E_8$), with its binary (ternary) Golay code construction.

It is widely known \cite{1989Koca}-\cite{Moxness2020-014} 
that the $E_8$ can be projected, mapped, or "folded" (as shown in Fig. \ref{fig:DynkinFold}) to two golden ratio $\varphi=\frac{1}{2}\left(1+\sqrt{5}\right)\approx1.618$ scaled copies of the 4 dimensional 120 vertex 720 edge $H_4$ 600-cell. Folding an 8D object into a 4D one can be done by projecting each vertex using its dot product with a 4$\times$8 matrix\cite{Moxness2014-006}. This produces $H_4$$\oplus$$\varphi H_4$, where $H_4$ is the binary icosahedral group 2$I$ of order 120, a subgroup of Spin(3). It covers $H_3$ as the full icosahedral group $I_h$ of order 120, a subgroup of SO(3). The binary icosahedral group is the double cover of the alternating group $A_5$.

Despite others'\cite{1998Koca}\cite{2016RSPSA.47250504D} recent attempts, the inverse morphism or "unfolding" from $H_4$ to $E_8$ is less trivial given that the matrix is not square and lacks an inverse. Yet, a real ($\mathbb{R}$) symmetric volume preserving Det($\mathbb{U}$)=1 rotation matrix(\ref{eqn:H4foldUnitary}) was derived in 2012 and documented\cite{Moxness2014-006}\cite{Moxness2018-011}\cite{Moxness2019-013}. The quadrant structure of $\mathbb{U}$ rotates $E_8$ into four 4D copies of $H_4$ 600-cells, with the original two (L)eft and (R)ight side unit scaled 4D copies related to the two L/R $\varphi$ scaled copies which we now identify as $H_{4}$(L$\oplus$R$\oplus$1$\oplus\varphi$). This traceless form of $\mathbb{U}$ has palindromic characteristic coefficients and provides for an explicit isomorphic mapping of $E_8$$\leftrightarrow$$H_{4}$(L$\oplus$R$\oplus$1$\oplus\varphi$). This involves using a bidirectional L$\leftrightarrow$R mapping function ($\mathtt{mapLR}$) and $\mathbb{U}^{-1}$(\ref{eqn:InvH4foldUnitary}). The process is described and visualized in Section \ref{sec:The palindromic unitary matrix}. It is interesting to note the exchange of 1$\leftrightarrow$$\varphi$ in $\mathbb{U}\leftrightarrow \mathbb{U}^{-1}$, excluding $-\varphi ^{2}$.

\begin{equation}
\label{eqn:H4foldUnitary}
\text{$ \mathbb{ U } $}\text{ = } \left(
\begin{array}{cccccccc}
1- \varphi& 0 & 0 & 0 & 0 & 0 & 0 & -\varphi ^{2} \\
0 & -1 & \varphi & 0 & 0 & \varphi & 1 & 0 \\
0 & \varphi & 0 & 1 & -1 & 0 & \varphi & 0 \\
0 & 0 & -1 & \varphi & \varphi & 1 & 0 & 0 \\
0 & 0 & 1 & \varphi & \varphi & -1 & 0 & 0 \\
0 & \varphi & 0 & 1 & -1 & 0 & \varphi & 0 \\
0 & 1 & \varphi & 0 & 0 & \varphi & -1 & 0 \\
 -\varphi ^{2} & 0 & 0 & 0 & 0 & 0 & 0 &1- \varphi \\
\end{array}
\right)/(2 \sqrt{\varphi})
\end{equation}

\begin{equation}
\label{eqn:InvH4foldUnitary}
\text{$\mathbb{U}^{-1}$}\text{=}\left(
\begin{array}{cccccccc}
\varphi-1& 0 & 0 & 0 & 0 & 0 & 0 & -\varphi ^{2} \\
0 & -\varphi & 1 & 0 & 0 & 1 & \varphi & 0 \\
0 & 1 & 0 & \varphi & -\varphi & 0 & 1 & 0 \\
0 & 0 & -\varphi &1 & 1 & \varphi & 0 & 0 \\
0 & 0 & \varphi & 1 &1& -\varphi & 0 & 0 \\
0 & 1 & 0 & \varphi & -\varphi & 0 & 1 & 0 \\
0 & \varphi & 1 & 0 & 0 & 1 & -\varphi & 0 \\
 -\varphi ^{2} & 0 & 0 & 0 & 0 & 0 & 0 &\varphi-1 \\
\end{array}
\right)/(2 \sqrt{\varphi})
\end{equation}

\subsection{Generating Polytopes}

\begin{figure}[!h]
\center
\includegraphics[width=200pt]{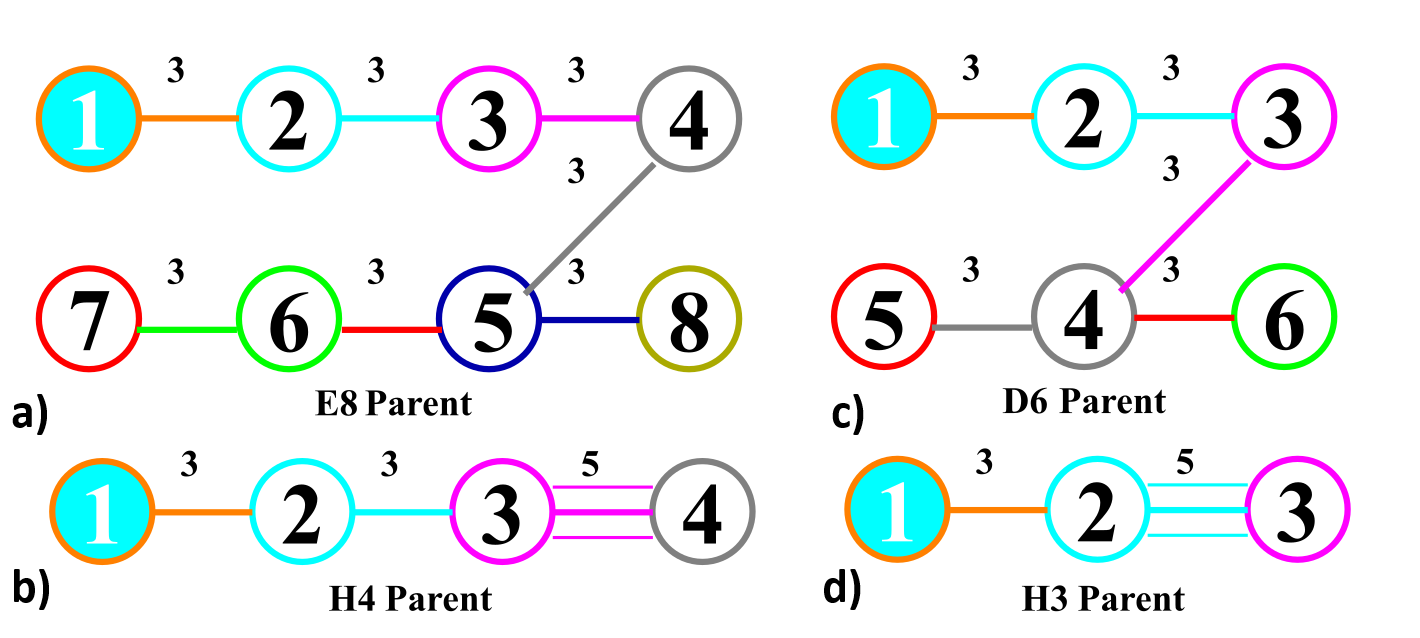}
\caption{\label{fig:DynkinFold}
a) $E_8$ Dynkin diagram in folding orientation\\
b) The associated Coxeter-Dynkin diagram of $H_4$\\
c) $D_6$ Dynkin diagram in folding orientation\\
d) The associated Coxeter-Dynkin diagram of $H_3$}
\end{figure}

The quaternion ($\mathbb{H}$) Weyl group orbit O($\Lambda$)=W($H_4$)=I of order 120 is constructed from the parent orbit (1000) of the Coxeter-Dynkin diagram for $H_4$ shown in Fig. \ref{fig:DynkinFold}b. This results in the 600-cell 4-polytope of order 120 labeled here and in \cite{Koca_2011} as I. In addition, $\mathbb{U}$ provides for a direct mapping from $E_8$ to four L$\oplus$R$\oplus$1$\oplus\varphi$ copies of the tri-rectified parent of $H_4$ (i.e. the filled node 1 is shifted right 3 times giving 0001), which is the 120-cell of order 600 labeled here and in \cite{Koca_2011} as J. Both of these 4-polytopes are shown in Appendix \ref{app:A} Figs. \ref{fig:421-concentric-hulls}-\ref{fig:J}. The detail of the quaternion Weyl orbit construction is described in Section \ref{sec:Quaternionic Weyl orbit construction}.

\begin{figure}[!ht]
\center
\includegraphics[width=255pt]{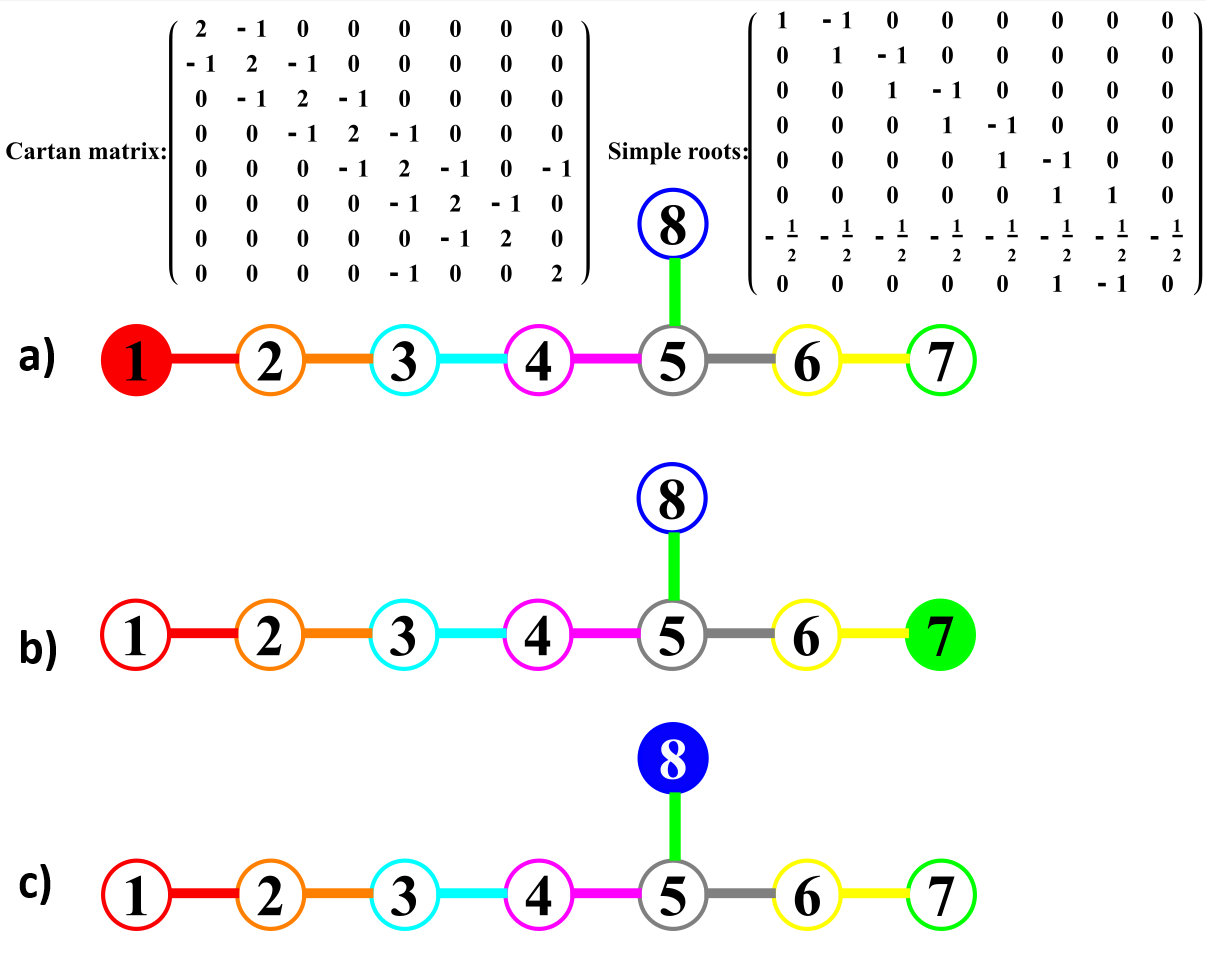}
\caption{\label{fig:421-dynkin} $E_8$ Dynkin diagrams a) $4_{21}$, b) $2_{41}$, c) $1_{42}$} Also shown are the Cartan and simple root matrices which correspond to the common Coxeter-Dynkin representation of the diagrams
\end{figure}

In addition to the 240 root $4_{21}$ $E_8$ 8-polytope identified by its Coxeter-Dynkin diagram in Fig. \ref{fig:421-dynkin}a, there are $2^8$ possible orbits using only 0's$\leftrightarrow$1's, empty$\leftrightarrow$filled, or ringed nodes of the $E_8$ Coxeter-Dynkin diagram, including the snub (00000000) orbit. Several other orbit permutations are commonly represented visually using the Petrie projection basis. They are the 2,160 root $2_{41}$ and 17,280 root $1_{42}$ 8-polytopes, which are constructed by generating the resulting roots by moving the filled (or ringed) node to each of the two other ends of the Dynkin diagram, as shown in Figs. \ref{fig:421-dynkin}b and \ref{fig:421-dynkin}c respectively.

\subsection{8D Platonic Rotation}

Interestingly from \cite{Moxness2019-013}, $\mathbb{U}$ can be generated using a combination of the unimodular matrices commonly used for Quantum Computing (QC) qubit logic, namely those of the 2 qubit CNOT (\ref{eqn:CNOT}) and SWAP (\ref{eqn:SWAP}) gates. Taking these patterns, combined with the recursive functions that build $\varphi$ from the Fibonacci sequence, it is straightforward to derive $\mathbb{U}$ from scaled QC logic gates.\cite{Moxness2020-014}

\begin{equation}
\label{eqn:CNOT}
\text{CNOT}\text{=}\left(
\begin{array}{cccc}
 1 & 0 & 0 & 0 \\
 0 & 1 & 0 & 0 \\
 0 & 0 & 0 & 1 \\
 0 & 0 & 1 & 0 \\
\end{array}
\right)
\end{equation}

\begin{equation}
\label{eqn:SWAP}
\text{SWAP}\text{=}\left(
\begin{array}{cccc}
 1 & 0 & 0 & 0 \\
 0 & 0 & 1 & 0 \\
 0 & 1 & 0 & 0 \\
 0 & 0 & 0 & 1 \\
\end{array}
\right)
\end{equation}

\subsection{2D and 3D Projection}

Projection of $E_8$ to 2D (or 3D) requires 2 (or 3) basis vectors $\{X,Y,Z\}$. For the Petrie projection shown in Fig. \ref{fig:E8Petrie}, we start with the basis vectors in (\ref{eqn:H4projVecs}), which are simply the two 2D Petrie projection basis vectors of the 600-cell (a.k.a. the Van Oss projection), with an optional 3rd (z) basis vector added for an interesting 3D projection\cite{Moxness2014-006}. 

\begin{equation}
\label{eqn:H4projVecs}
\begin{array}{cccccc}
\text{x=}\{&0, &\varphi 2\text{Sin}\frac{2\pi}{15},&2\text{Sin}\frac{2\pi}{15},&0,&0, 0, 0, 0\} \\
\text{y=}\{&-\varphi 2\text{Sin}\frac{2\pi}{30}, &0, &0, &1,&0, 0, 0, 0\}\\
\text{z=}\{&1, &0, &0, &\varphi 2\text{Sin}\frac{2\pi}{30},&0, 0, 0, 0\}\\
\end{array}
\end{equation}

 $\{X,Y,Z\}=\mathbb{U}.\{x,y,z\}$ as shown in (\ref{eqn:E8projVecs}).
\begin{equation}
\label{eqn:E8projVecs}
\begin {array} {cccccccc} 
\text{X=}\{0 & .252 & .427 & -.319 & .319 & .427 & .781 & 0\} \\ 
\text{Y=}\{.0821 & 0 & - .393 & .636 & .636 & .393 & 0 & .348\} \\
\text{Z=}\{-.242 & 0 & -.132 & .215 & .215 & .132 & 0 & -1.03\} \\
\end {array}
\end{equation}

\subsection{3D Platonic Solid Projection}

\begin{figure}[!ht]
\center
\includegraphics[width=100pt]{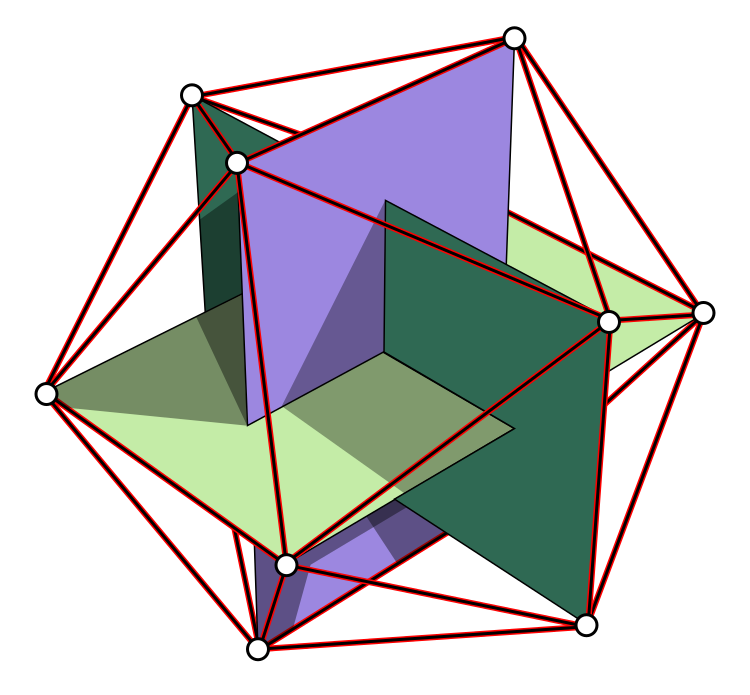}
\caption{\label{fig:Icosahedron}The icosahedron formed from 3 mutually-perpendicular golden rectangles}
\end{figure}

This basis is derived from the icosahedral symmetry of the $H_3$-based Platonic solid. The twelve vertices of the icosahedron can be decomposed into three mutually-perpendicular golden rectangles (as shown in Fig. \ref{fig:Icosahedron}), whose boundaries are linked in the pattern of the Borromean rings. Rows (or columns) 2-4 (or 5-8) of $\mathbb{U}$ contain 6 of the 12 vertices of this icosahedron, including 2 at the origin with the other 6 of 12 icosahedron vertices being the antipodal reflection of these through the origin. These 2 (or 3) rows can then used as a kind of ``Platonic solid projection prism'' to form the 2 (or 3) 8D basis vectors used in the 2D (or 3D) projection of $4_{21}$, $2_{41}$, and $1_{42}$. 

\begin{figure}
\center
c)\includegraphics[width=180pt]{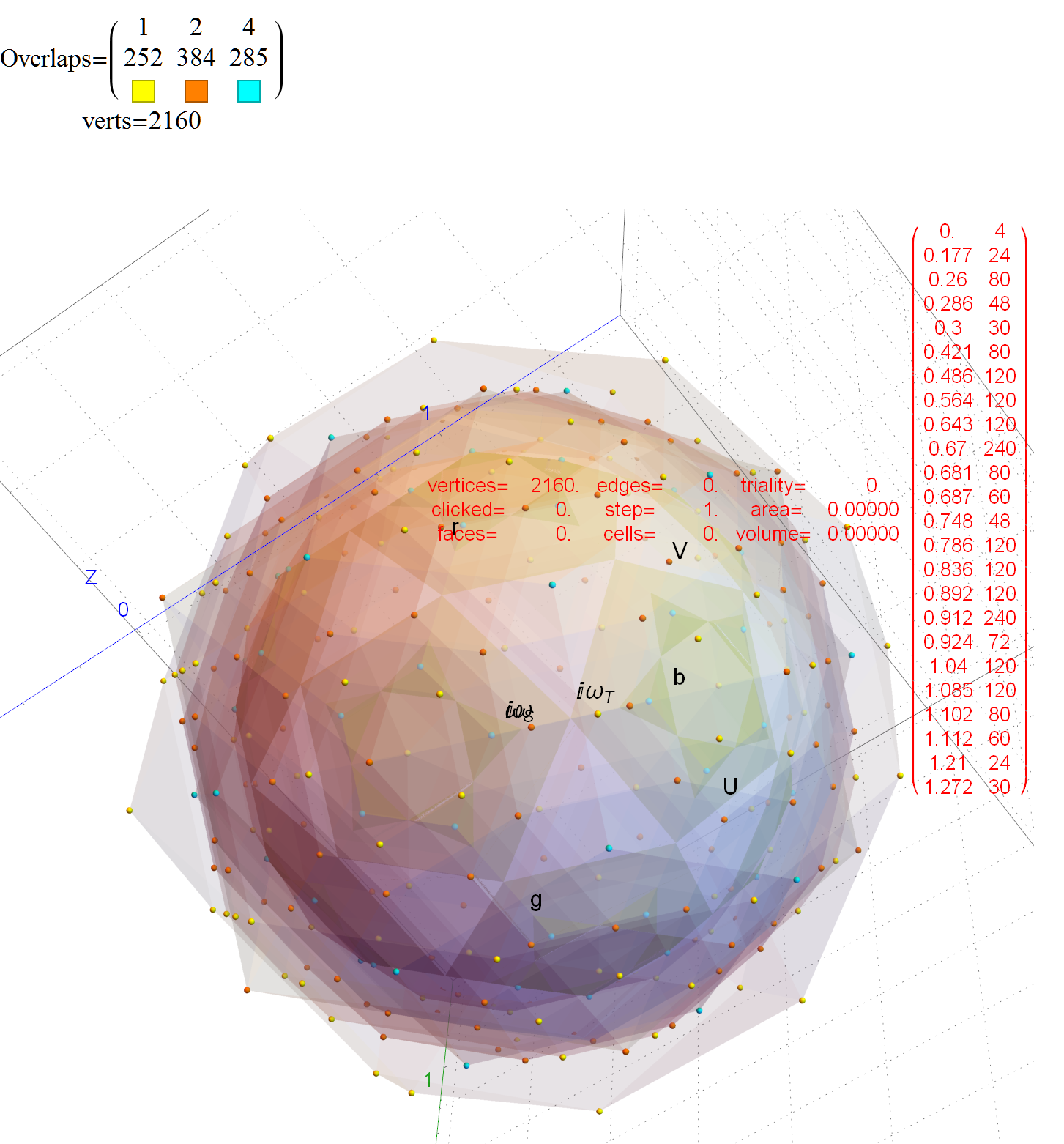}
\caption{\label{fig:E8-241}
$2_{41}$ projections of its 2,160 vertices\\
a) 2D to the $E_8$ Petrie projection using basis vectors X and Y from (\ref{eqn:E8projVecs}) with 8-polytope radius 2$\sqrt{2}$ and 69,120 edges of length $\sqrt{2}$.\\
b) 3D projections with vertices sorted and tallied by their 3D norm generating the increasingly transparent hulls for each set of tallied norms. Notice the last two outer hulls are a combination of two overlapped Icosahedrons (24) and a Icosidodecahedron (30).\\
c) Combined 3D hulls with the overlapping vertices color coded by overlap count. Also shown is a list (in red) of the normed hull distance and the number of vertices in the group.}
\end{figure}

\begin{figure}
\center
a)\includegraphics[width=185pt]{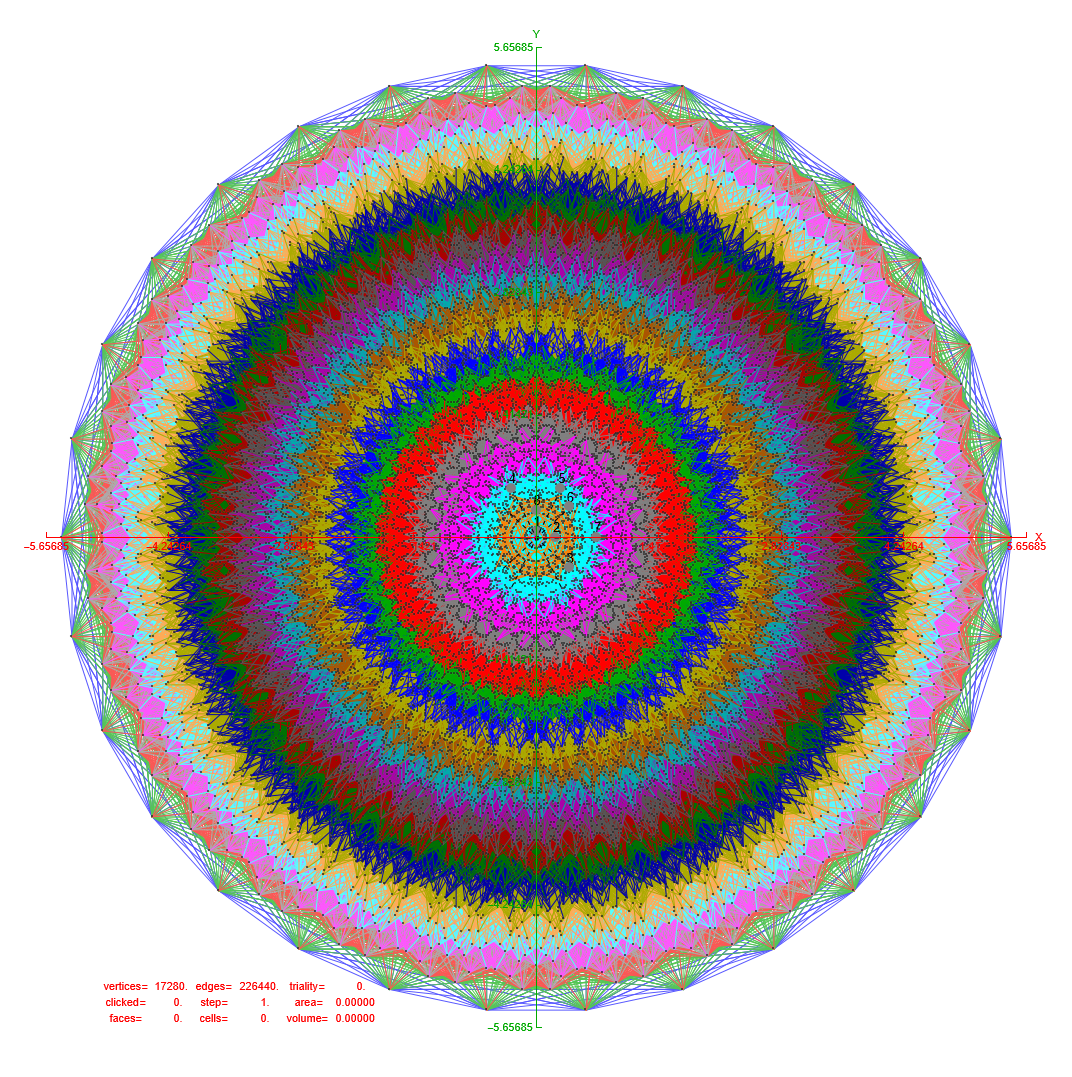}

b)\includegraphics[width=180pt]{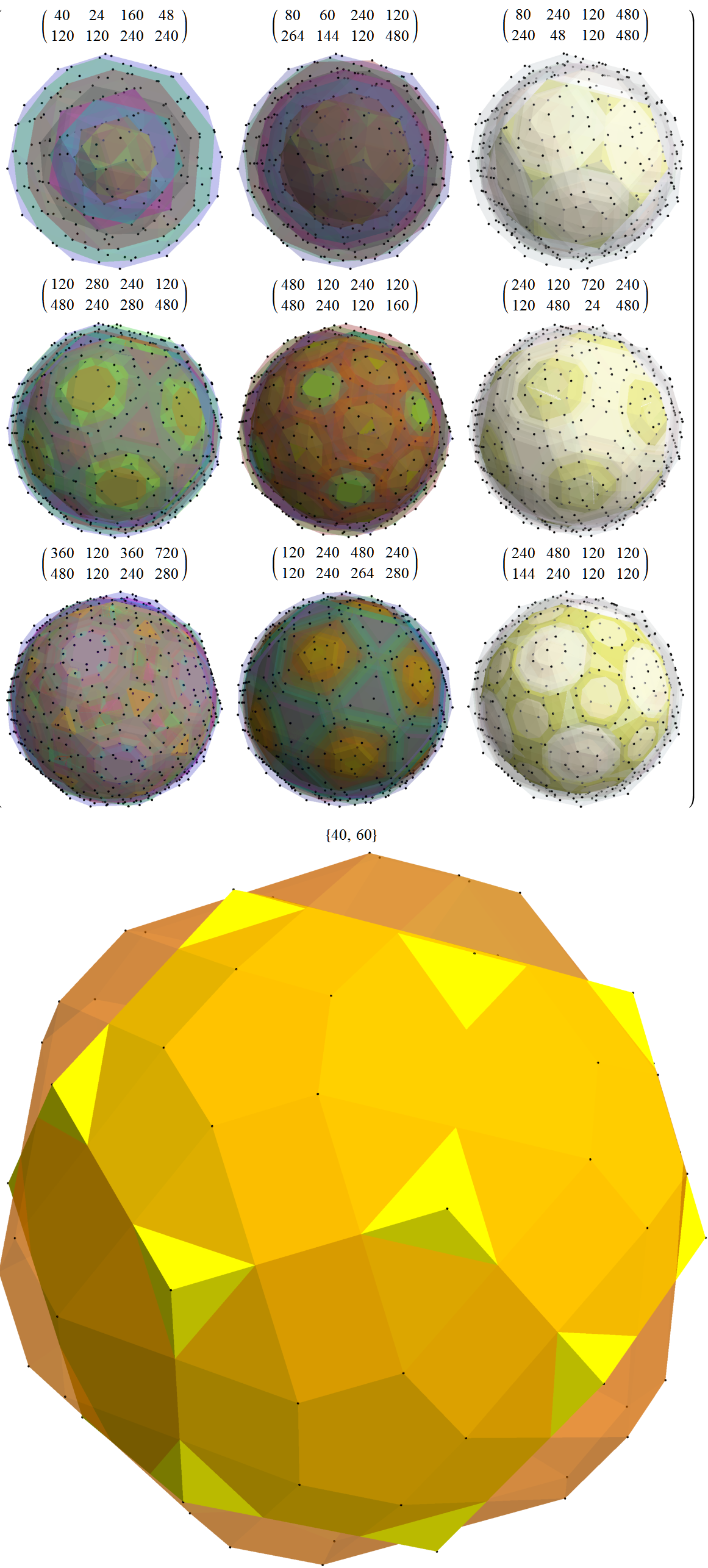}

\caption{\label{fig:E8-142}$1_{42}$ projections of its 17,280 vertices\\
a) 2D to the $E_8$ Petrie projection using basis vectors X and Y from (\ref{eqn:E8projVecs}) with 8-polytope radius 4$\sqrt{2}$ and 483,840 edges of length $\sqrt{2}$ (with 53$\%$ of inner edges culled for display clarity).\\
b) 3D projections with vertices sorted and tallied by their 3D norm generating the increasingly transparent hulls for each set of tallied norms. Notice the last two outer hulls are a combination of two overlapped Dodecahedra (40) and a irregular Rhombicosidodecahedron (60).}
\end{figure}

Orthogonal projection to 3D after $\mathbb{U}$ folding (i.e. selecting one of 56 unique subsets of any 3 dimensions, here we use $\{1,2,3\}$) manifests a large number of concentric hulls with Platonic and Archimedean solid related structures. The eight projected 3D hulls of $4_{21}$ include two $\varphi$ scaled sets of four hulls from two 600-cells ($H_4 \oplus \varphi H_4$) as shown in Appendix \ref{app:A} Fig. \ref{fig:421-concentric-hulls}. $2_{41}$ and $1_{42}$ projections of $E_8$ are shown in Figs. \ref{fig:E8-241}-\ref{fig:E8-142}. 

\begin{figure}
\center
a)
\includegraphics[width=120pt]{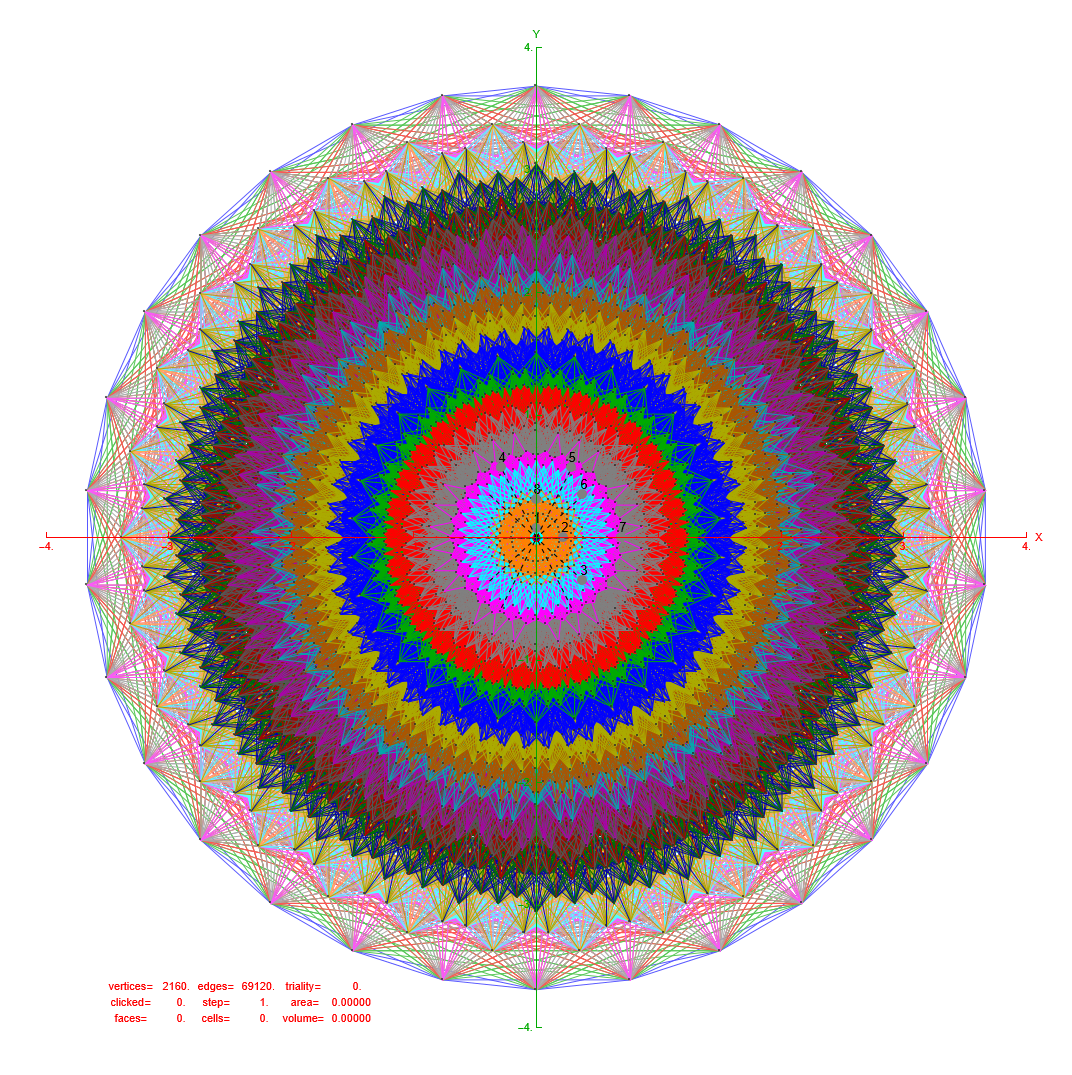}\\
\end{figure}

\begin{figure}
\center
b)\includegraphics[width=190pt]{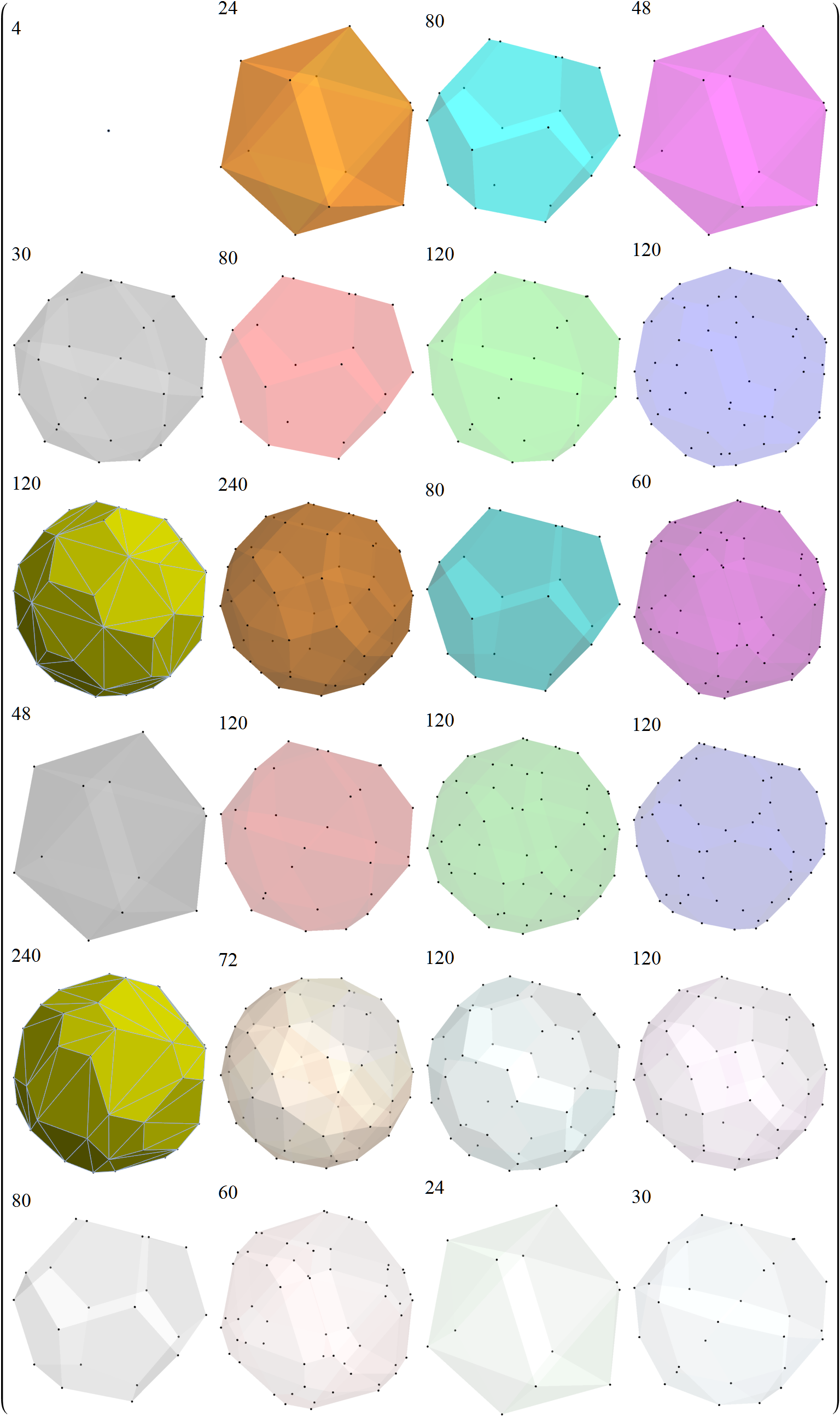}\\
\end{figure}

\section{The palindromic unitary matrix}
\label{sec:The palindromic unitary matrix}

The particular maximal embedding of $E_8$ at height 248 that we are interested in for this work is shown in Appendix \ref{app:C} Fig. \ref{fig:CombinedEmbeddings} as the special orthogonal group of SO(16)=$D_8$ at height (120=112+4+4)+128’, where 112 is interpreted as the subgroup embeddings of SO(8)$\otimes$SO(8)=$D_4$$ \otimes$$D_4$ and 128' is interpreted as symplectic subgroup embeddings of $C_8$ where Sp(8)$\otimes$Sp(8)=$C_4$$\otimes$$C_4$ at height 136=128+4+4. These selected embeddings correspond to the 112 integer $D_8$ vertices and the 128 half-integer $BC_8$ vertices given by SRE $E_8$, in addition to the 8$\oplus \overline{8}$ generator roots for a total of $2^8$. This is in 1::1 correspondence with the canonical root vertex ordering from the 9th row of the palindromic Pascal triangle $\{1,8,28,56,35\overline{35},\overline{56},\overline{28},\overline{8},\overline{1}\}$, where each entry in the list gives the number of vertices that alternate between half-integer $BC_8$ and integer $D_8$ vertex sets, with the right 5 overbar sets of 128 vertices being the negated vertices of the left 5 sets of 128 in reverse order.

These embeddings have an isomorphic connection to $\mathbb{U}$ and provide the $E_8$$\leftrightarrow$$H_{4}$(L$\oplus$R$\oplus$1$\oplus\varphi$) mapping via $\mathtt{mapLR}$. The \textit{Mathematica}\textsuperscript{TM} code for $\mathtt{mapLR}$ and the code to validate the $E_8\leftrightarrow H_4$ isomorphism is shown in Appendix \ref{app:D} Fig. \ref{fig:Isomorphism-code}. It demonstrates that $E_8$ rotates into four 4D copies of $H_4$ 600-cells, with the original two (L)eft side $\varphi$ scaled 4D copies related to the two (R)ight side unscaled 4D copies. 
testtest
Due to the palindromic structure of $\mathbb{U}$, the $H_{4L}$ and $H_{4R}$ are also palindromic with each R vertex being the reverse order of the L vertex, along with $\mathtt{mapLR}$ exchanges in the (S)nub 24-cell vertices. For each L vertex that is not a member of the (T)etrahedral group's self-dual $D_4$ 24-cell (or $\varphi$T), the R vertex will be a member of the scaled $\varphi$S (or S) respectively. This is due to the exchange of $\varphi^{3/2}$$\leftrightarrow$$\varphi^{-3/2}$ in $\mathtt{mapLR}$ which changes the norm (i.e. to/from a small norm=$1/\sqrt{\varphi}$ or a large norm=$\sqrt{\varphi}$). The 24-cell T vertices are unaffected by $\mathtt{mapLR}$ exchange and have L and R vertex values of the same norm and palindromic opposite entries, with the larger $\varphi H_4$ having the same signs and the smaller unit scaled $H_4$ having opposite signs.

\begin{figure}[!ht]
\center
\includegraphics[width=265pt]{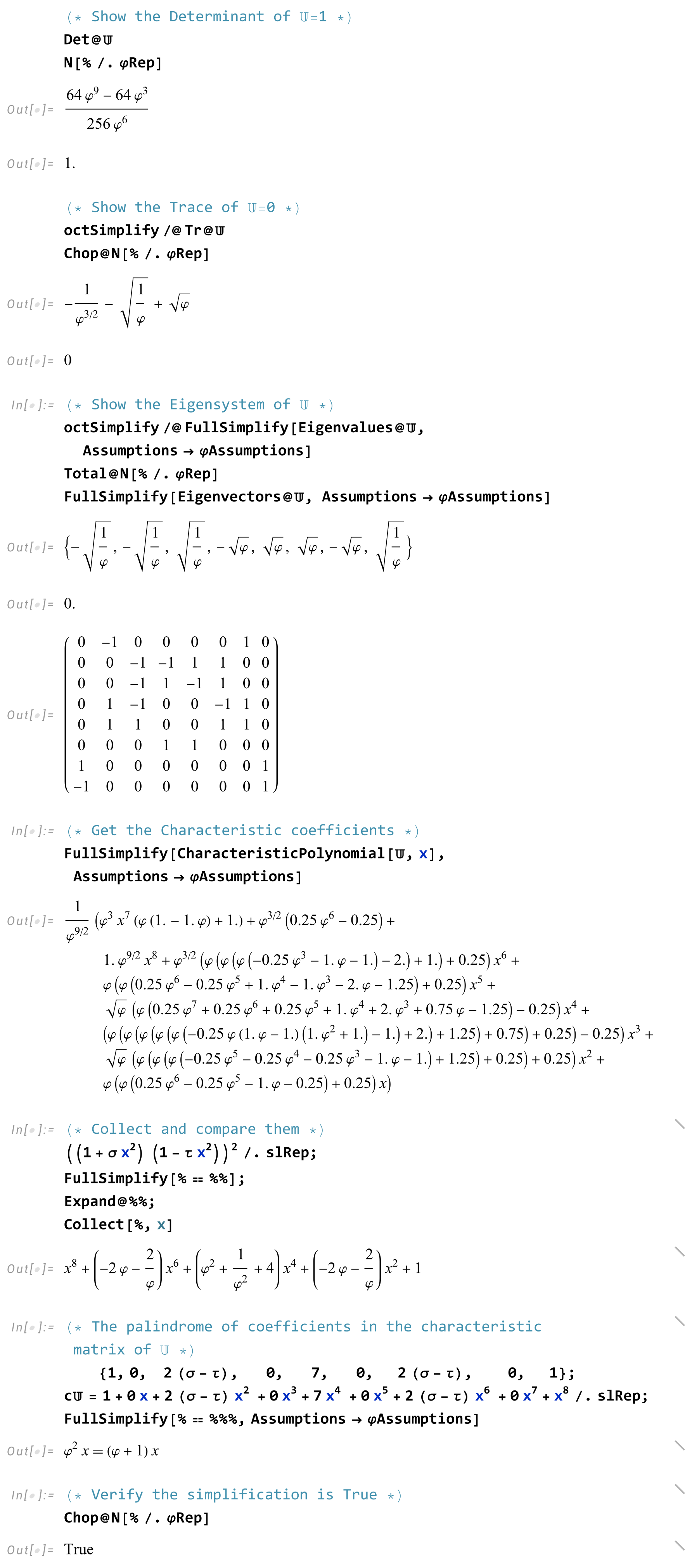}
\caption{\label{fig:U-Characteristic-Coefficients}The trace, determinant, Eigenvalues, Eigenvector matrix, and characteristic polynomial coefficients of $\mathbb{U}$} 
\end{figure}

\begin{figure}[!ht]
\center
\includegraphics[width=280pt]{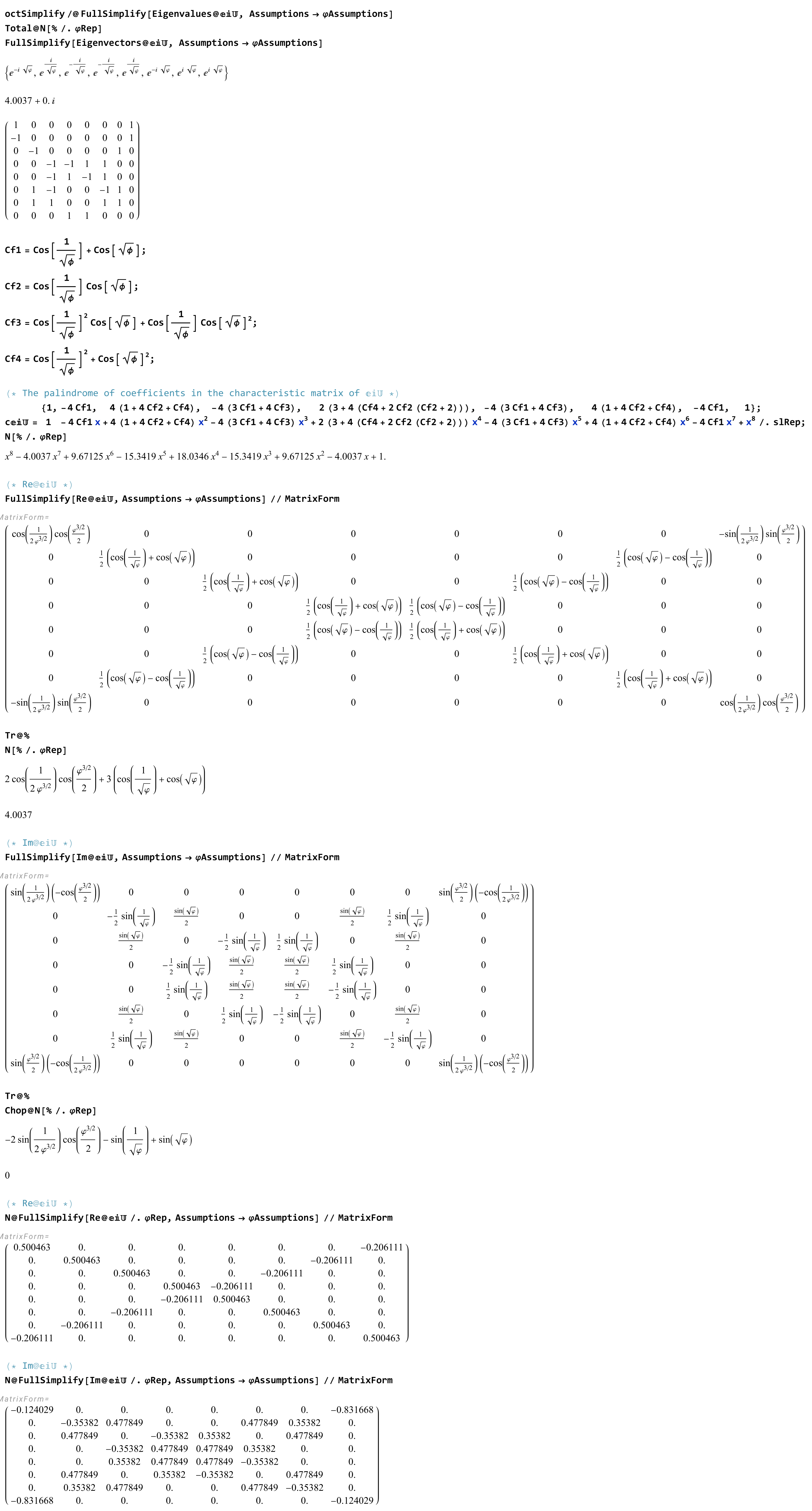}
\caption{\label{fig:eIU-Characteristic-Coefficients}The Eigenvalues, Eigenvector matrix, and characteristic polynomial coefficients of the unitary form of $\mathbb{U}$ as $e^{\text {i$\mathbb{U}$}}$ showing a Tr@Re@$e^{\text {i$\mathbb{U}$}}\approx 4$ and a traceless imaginary part}
\ \\
\ \\
\ \\
\ \\
\end{figure}

It is clear that $\mathbb{U}$ is traceless, but it is not unitary. Since $\mathbb{U}$ is Hermitian, it is easily made unitary as $e^{\text {i$\mathbb{U}$}}$. While that is unitary it is not traceless, so it is not an $A_7$ group SU(8) symmetry. For the identification of their palindromic characteristic polynomial coefficients, see Figs. \ref{fig:U-Characteristic-Coefficients}-\ref{fig:eIU-Characteristic-Coefficients}. 

See Appendix \ref{app:D} Figs. \ref{fig:Isomorphism-output-H4}-\ref{fig:Isomorphism-output-H4phi} showing the detail of the $E_8$$\leftrightarrow$$H_{4}$(L$\oplus$R$\oplus$1$\oplus\varphi$) isomorphism and the patterns within their respective vertex roots.
\ \\
\ \\
\ \\
\ \\

\section{Quaternionic Weyl orbit construction}
\label{sec:Quaternionic Weyl orbit construction}

\begin{figure}[!ht]
\center
\includegraphics[width=160pt]{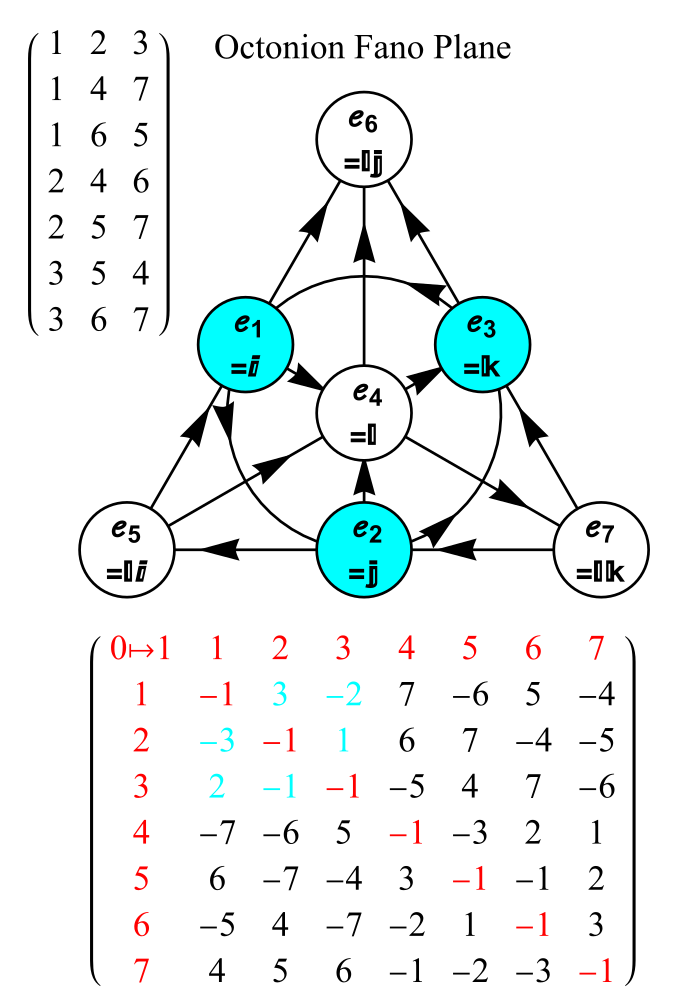}
\caption{\label{fig:SetUpFano}The selected octonion Fano plane mnemonic and multiplication table based on its 7 structure constant triads . The first triad (123) defines standard convention for quaternions.}
\end{figure}

\begin{figure}[!ht]
\center
\includegraphics[width=160pt]{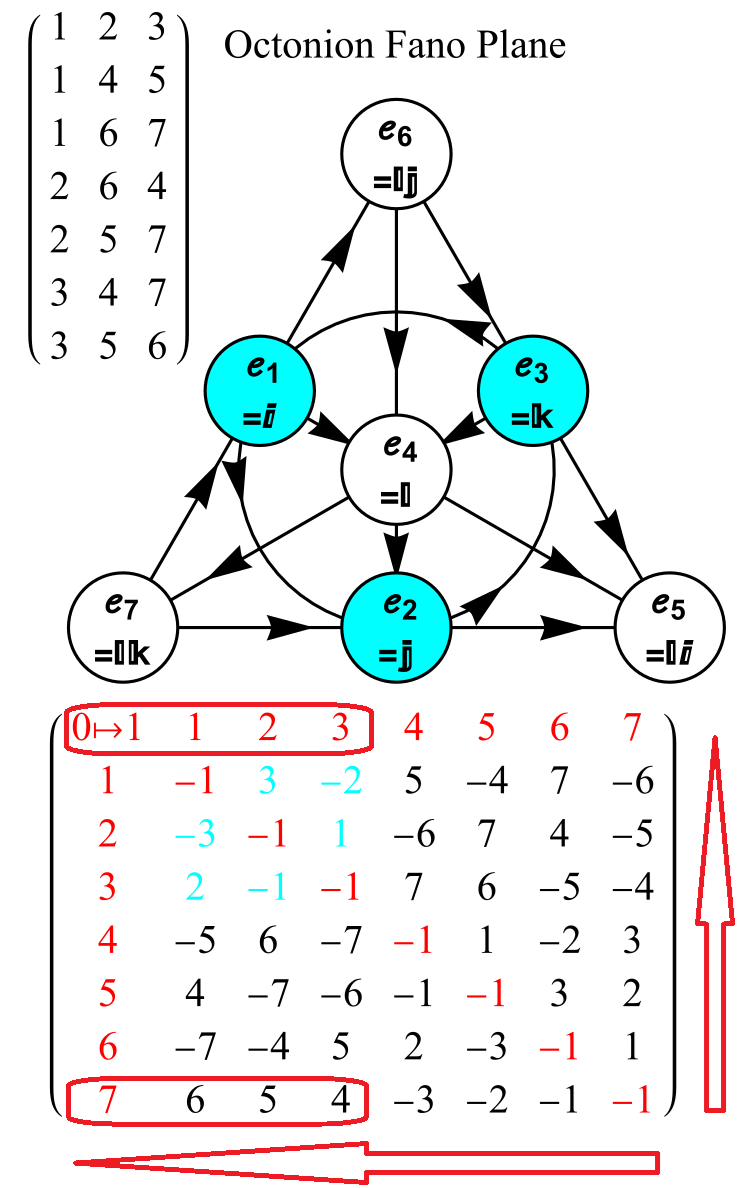}
\caption{\label{fig:palindromic-Fano}An alternative set of structure constant triads, octonion Fano plane mnemonic, and multiplication table, with decorations showing the palindromic multiplication.}
\end{figure}

The content within this paper was generated using a computational environment the author has written in \textit{Mathematica\textsuperscript{TM} by Wolfram Research, Inc.}. In order to deal effectively with quaternions, it supplants the native Quaternion package with a more flexible symbolic octonion ($\mathbb{O}$) capability. This allows for the selection of a multiplication table from any of the 480 possible octonion tables, including their split and bi-octonion forms. It also handles the sedenion forms as well and has been used to verify the octonion forms of $E_8$ from Koca\cite{1989Koca}, Dixon\cite{dixon2010integral}, Pushpa and Bisht\cite{Pushpa_2012}, R. A. Wilson, Dray, and Monague\cite{Wilson_2023}, including the complexified octonions of Günaydin-Gürsey\cite{10.1063/1.1666240} and Furey\cite{Furey_2018}. To ensure that our quaternion (and bi-quaternion) math is consistent with the standard multiplication convention related to quaternions, we need to select one of the 48 octonions with a first triad of 123 and a Cayley-Dickson construction where $e_4$-$e_7$ quadrant multiplication remains within the quadrant. See Fig. \ref{fig:SetUpFano} showing the selected triads, Fano plane, and multiplication table of the octonion used in this and several of the referenced papers\footnote[1]{It is interesting to note that this particular octonion is close to (but not) palindromic. Using an algorithmic identification and construction of all of the possible 480 unique permutations of octonions\cite{Moxness2013-005}, we find that a small change in triads to $\{$123,145,167,264,257,347,356$\}$ with 5$\leftrightarrow$7 ordering swaps creates a palindromic $E_8$. This octonion is shown in Fig. \ref{fig:palindromic-Fano}.}.

\begin{figure}[!ht]
\center
\includegraphics[width=260pt]{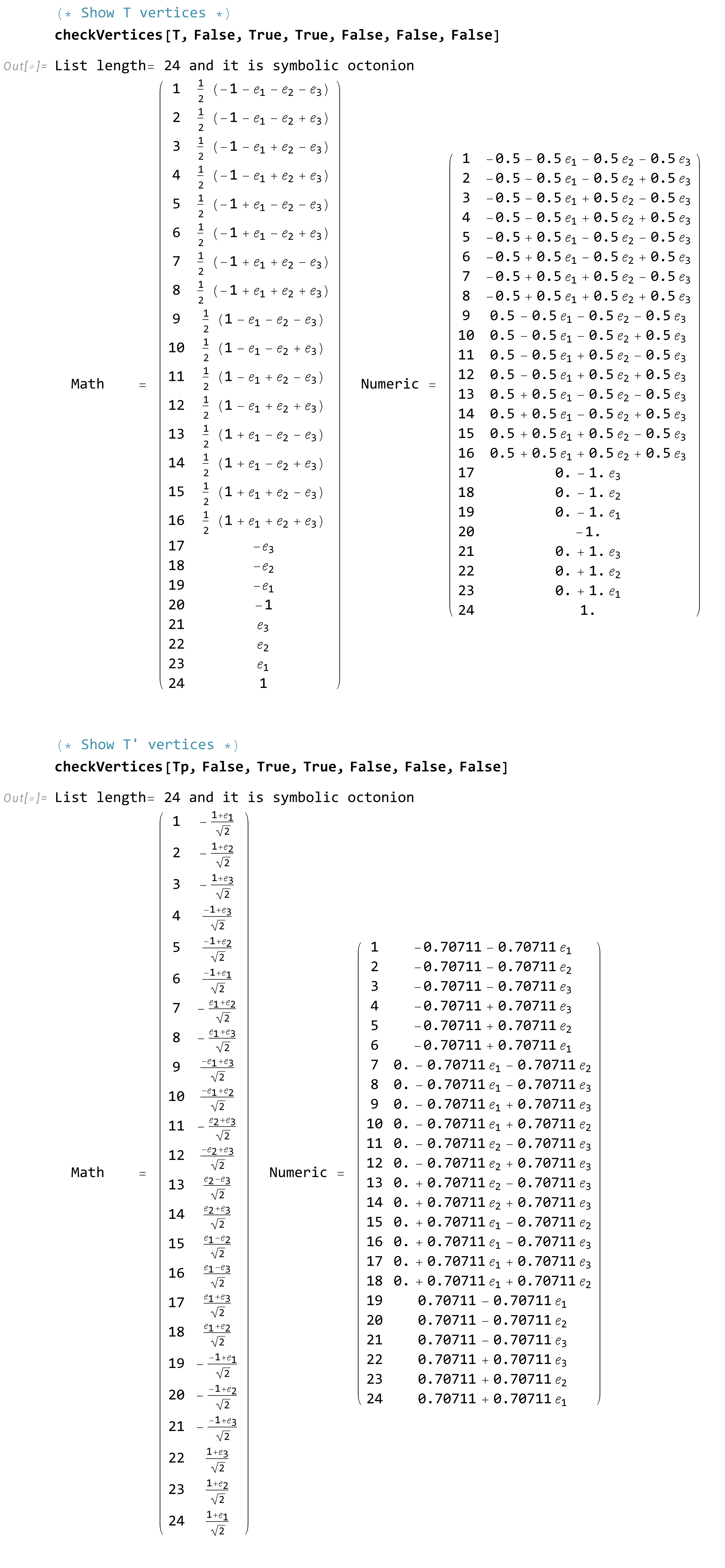}
\caption{\label{fig:TandTp}The values of the $D_4$ 24-cell T and its alternate T'\\
}
\ \\
\ \\
\ \\
\ \\
\ \\
\end{figure}

\begin{figure}[!ht]
\center
\includegraphics[width=220pt]{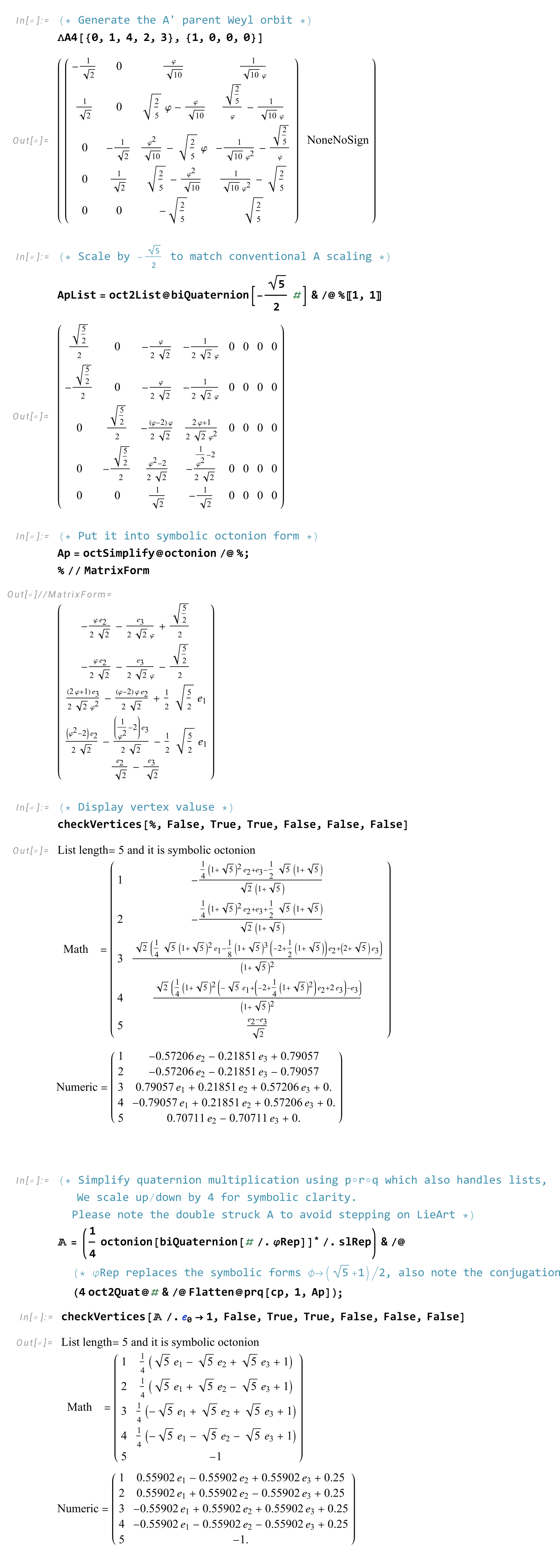}
\caption{\label{fig:AfromAp} Explicit \textit{Mathematica\textsuperscript{TM}} computation of A from the $\mathtt{\Lambda A4[\Lambda\_,orbit\_]}$ generated A'}
\end{figure}

\begin{figure}[!ht]
\center
\includegraphics[width=200pt]{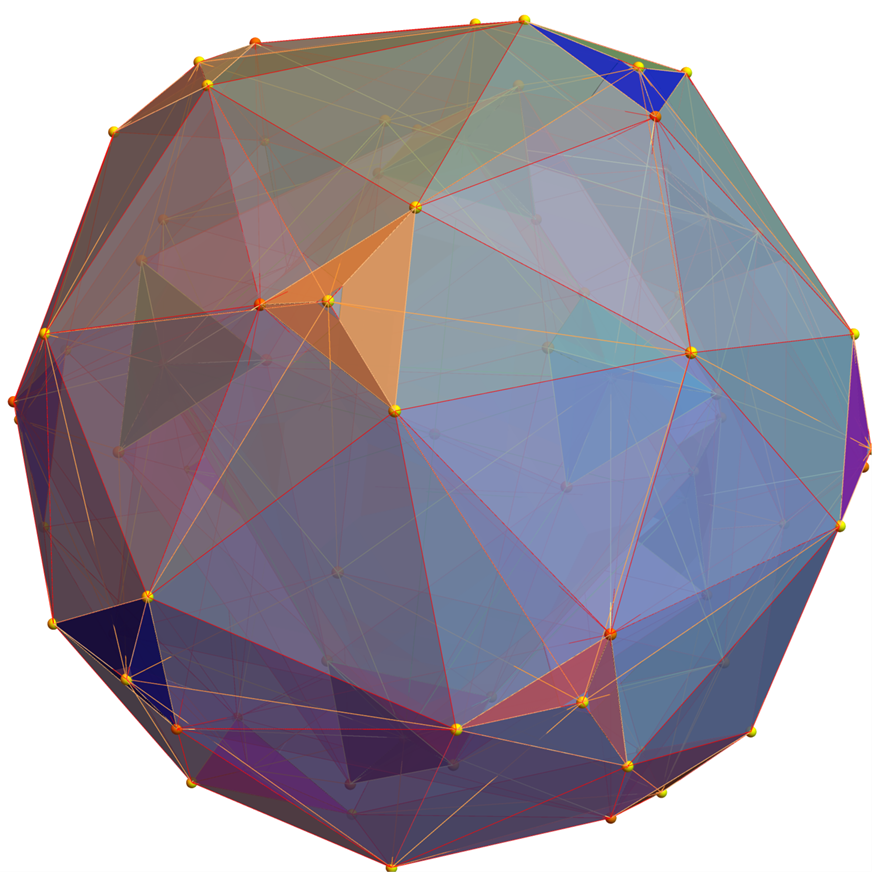}
\caption{\label{fig:Dual_Snub_24_Cell}Visualization of the 144 root vertices of S'+T+T' now identified as the dual snub 24-cell}
\end{figure}

It has been shown that the 3D symmetry groups of $A_3$, $B_3$, and $H_3$\cite{Koca_2011} and 4D symmetry groups of $A_4$, $D_4$, $F_4$, and $H_4$ are related to the higher dimensional groups of $D_6$ and $E_8$\cite{Koca_2012}\cite{2016RSPSA.47250504D}. A quaternionic Weyl group orbit O($\Lambda$)=W($H_4$)=I of order 120 can be constructed from $H_3$ which generates some of the Platonic, Archimedean and dual Catalan solids shown in Appendix \ref{app:B} Fig. \ref{fig:Archimedean-and-Catalan-solids}, including their irregular and chiral forms\cite{koca2011chiral}. The polytopes for a particular orbit of O($\mathtt{\Lambda}$)=W($\mathtt{group}$) are generated using a function $\mathtt{\Lambda[group\_,orbit\_,perm\_:"Rotate"]}$, where $\mathtt{perm}$ can be one of 18 combinations of sign and position permutation functions (e.g. "oSign" gives all odd sign permutations and cyclic rotations of position and the default "Rotate" gives all sign permutations of cyclically rotated positions). The first column in these figures show the set of calls to the $\mathtt{\Lambda}$ function. This same method is used to generate the $H_4$-based 4-polytopes of the 120-cell and 600-cell shown in Appendix \ref{app:A} Figs. \ref{fig:421-concentric-hulls}-\ref{fig:J}.

The $A_3$ in $A_4$ group embedding of SU(5)$\supset$SU(4)$\otimes$$U_1$\cite{Koca_2012} are shown in Appendix \ref{app:C} Fig. \ref{fig:A4-A3-SU5xSU4xU1-1024-Koca} in combination with these 3 and 4-polytope visualizations.
\footnote[2]{In the methods and coding descriptions, since Mamone\cite{Mamone-sym2031423} identifies the 5-cell as S, but Koca uses S to identify the (S)nub 24-cell (a convention which we use here), Mamone's $A_4$-based 5-cell is now identified as A which is the 4D version of the tetrahedron.}

We identify the rectified parent orbit (0100) of W($D_4$) as the self-dual 24-cell T, which is the combination of the 4D octahedron (aka. 16-cell) and the 4D cube (aka. 8-cell with a 3D hull of the cuboctahedron derived from the tri-rectified (0001) W($BC_4$)). Due to the W($D_4$) Coxeter-Dynkin diagram triality symmetry, T' is identified with any of 3 end nodes as parent and others as bi-rectified and tri-rectified orbits $\{$(1000), (0010), (0001)$\}$ each with 8 vertices of 2-component (vector) quaternions and has a 3D hull of the rhombic dodecahedron. See Fig. \ref{fig:TandTp} for their specific symbolic and numeric values. Of course, it has also been shown that the root system of $F_4=T\oplus$T'.

From T (and T') we can take any one vertex to define a c (and c'=cp) respectively. For this paper, we use as an example c=$t_1$ from eq. (18) from Koca\cite{Koca_2011} T (and T') shown as \#13 in Fig. \ref{fig:TandTp} such that c=$\frac{1}{2} \left(1+\mathit{e}_1-\mathit{e}_2-\mathit{e}_3\right)$ (and c'=$\frac{\mathit{e}_2-\mathit{e}_3}{\sqrt{2}}$). Here c' is used with A' to generate the parent W($A_4$), or simply A as the 5-cell\cite{Koca_2011}. Specifically, A=$(c'\circ A')^{*}$ with A'=$\mathtt{\Lambda A4[\{0,1,4,2,3\},\{1,0,0,0\}]}$.
\footnote[3]{The 4-polytopes for a particular orbit of O($\mathtt{\Lambda}$)=W($\mathtt{group}$) are generated using a function $\mathtt{\Lambda[group\_,orbit\_,perm\_]}$ which is called by $\mathtt{\Lambda A4[\Lambda\_,orbit\_]}$ for the subgroup embeddings in $A_4$ as described in \cite{Koca_2012}. In addition, $\mathtt{SmallCircle}$ ($\circ$) is the symbolic operator for quaternion (octonion) multiplication that operates across lists, along with the expected symbolic exponentials (* and $\dag$) for Conjugate and ConjugateTranspose respectively. The function $\mathtt{prq[p\_,r\_,q\_,left_:False]:=If[left,(p\circ r)\circ q,p\circ (r\circ q)]}$ implements the operation of [p,q]:r from eq. (6) in \cite{Koca_2011}, which is defined for any combinations of inputs as elements or lists in order to add flexibility to quaternion and octonion operators, including left or right (default) non-commutative multiplication ordering. Other operators are also available for scalar product+($\oplus$), scalar product-($\ominus$), commutator($\odot$), anti-commutator($\wedge$), derivation($\Box$), Kronecker product($\otimes$), and $\mathtt{octExp}$ for exponential powers of octonions.} See Fig. \ref{fig:AfromAp} for the explicit \textit{Mathematica\textsuperscript{TM}} computation related to A and A'.

The snub orbit (0000) of W($D_4$) will generate the vertices of the snub 24-cell or S=I-T, as with the alternate snub 24-cell S'=I'-T' as shown in (\ref{eqn:S}) and (\ref{eqn:Sp}). We can generate S (or S') by taking the odd (or even) sign and cyclic position permutations of a seed quaternion p$\in$S (or S') to be assigned to $\alpha$ (or $\beta$) for generating S (or S') respectively. There are only 48 that satisfy the necessary constraint where a unit normed $p^5=\pm 1$. Those quaternions that satisfy the constraint are identified with an * in Appendix \ref{app:D}. For this paper, we selected from the 96 permutations of S $\alpha=\frac{1}{2} \left(\frac{1}{\varphi }+\varphi \mathit{e}_2+\mathit{e}_1\right)$ (and S' for $\beta=\frac{-\varphi -\frac{\mathit{e}_2}{\varphi }+\sqrt{5} \mathit{e}_1}{\sqrt{8}}$). This process of generating the snub 24-cell can be visualized as generating four quaternion 4D rotations of T (and T'). The 3D hulls of I'are shown in Fig. \ref{fig:Ip}. 

\begin{equation}
\label{eqn:S}
\begin{array}{l}
S=I-T=\sum _{i=1}^4 \alpha ^i\circ T\\
\mathtt{or}\\
I=\mathtt{prq[}\alpha^{0-4}\mathtt{,1,T]}
\end{array}
\end{equation}

\begin{equation}
\label{eqn:Sp}
\begin{array}{l}
S'=I'-T'=\sum _{i=1}^4 \beta ^i\circ T'\\
\mathtt{or}\\
I'=\mathtt{prq[}\beta^{0-4}\mathtt{,1,T']}
\end{array}
\end{equation}

The 3D hulls for one copy of I (or $\varphi$I) are represented in Fig. \ref{fig:421-concentric-hulls} hulls $\{$2,3,5$\}$ (or $\{$6,7,8$\}$) respectively plus 1/2 of the vertices in hull 4. The vertex values of I are listed in either of the center columns of Appendix \ref{app:D} Fig. \ref{fig:Isomorphism-output-H4} or Fig. \ref{fig:Isomorphism-output-H4phi}.

Koca\cite{Koca_2011} has also identified the dual to the snub 24-cell as being made up of the 144 root vertices of S'+T+T'. This 4-polytope is visualized in Fig. \ref{fig:Dual_Snub_24_Cell}.

The equations for the generation of J (and J') are shown in (\ref{eqn:J}) and (\ref{eqn:Jp}). As it was for I (and I') vertices each mapping to 5 quaternion rotations of T (and T'), J (and J') vertices each map to 5 quaternion rotations of I (and I') or 25 quaternion rotations of T (and T'). Given the isomorphism between each $E_8$ root vertex and 4 copies of I (i.e. L and R each at unit and $\varphi$ scales) as demonstrated in Section \ref{sec:The palindromic unitary matrix}, this means quaternionic Weyl orbit construction, when used with $\mathbb{U}$ and $\mathtt{mapLR}$, provides for an explicit map between each of the 240 $E_8$ root vertices and 10 J (or J') vertices (i.e. 10=2(L$\oplus$R)$\times$5 quaternion rotations of each I (or I') vertex).

\begin{equation}
\label{eqn:J}
\begin{array}{l}
J=\sum _{i=0}^4 c'\circ \bar{\alpha }^{\text{$\dagger $i}} \circ \alpha ^i\circ T\\
\mathtt{or}\\
J=\mathtt{prq[A',}\alpha^{0-4}\mathtt{,T]}
\end{array}
\end{equation}

\begin{equation}
\label{eqn:Jp}
\begin{array}{l}
J'=\sum _{i=0}^4 c\circ \bar{\beta }^{\text{$\dagger $i}} \circ \beta ^i\circ T'\\
\mathtt{or}\\
J'=\mathtt{prq[A',}\beta^{0-4}\mathtt{,T']}
\end{array}
\end{equation}

See Figs. \ref{fig:J}-\ref{fig:Jp} for the 120-cell (J) and its alternate (J') as generated by J=$\mathtt{prq[A',1,I]}$ and J'=$\mathtt{prq[A',1,I']}$ respectively.
\ \\
\ \\

\section{Conclusion}

This paper has given an explicit isomorphic mapping from the 240 $\mathbb{R}^{8}$ root $E_8$ Gosset $4_{21}$ 8-polytope to two $ \varphi$ scaled copies of the 120 root $H_4$ 600-cell quaternion 4-polytope using $\mathbb{U}$. It has also shown the inverse map from a single $H_4$ 600-cell to $E_8$ using a 4D$\hookrightarrow$8D chiral L$\leftrightarrow$R mapping function, $\varphi$ scaling, and $\mathbb{U}^{-1}$. This approach has shown that there are actually four copies of each 600-cell living within $E_8$ in the form of chiral $H_{4L}$$\oplus$$\varphi H_{4L}$$\oplus$$H_{4R}$$\oplus$$\varphi H_{4R}$ roots. In addition, it has demonstrated a quaternion Weyl orbit construction of $H_4$-based 4-polytopes that provides an explicit map from $E_8$ to four copies of the tri-rectified Coxeter-Dynkin diagram of $H_4$, namely the 120-cell of order 600. Taking advantage of this property promises to open the door to as yet unexplored chiral $E_8$-based Grand Unified Theories or GUTs. It is anticipated that these visualizations and connections will be useful in discovering new insights into unifying the mathematical symmetries as they relate to unification in theoretical physics. 

\begin{acknowledgments}

I would like to thank my wife for her love and patience and those in academia who have taken the time to review this work.

\end{acknowledgments}

\bibliography{The_Isomorphism_of_H4_and_E8}

\appendix

\begin{figure}[!ht]
\center
\includegraphics[width=525pt]{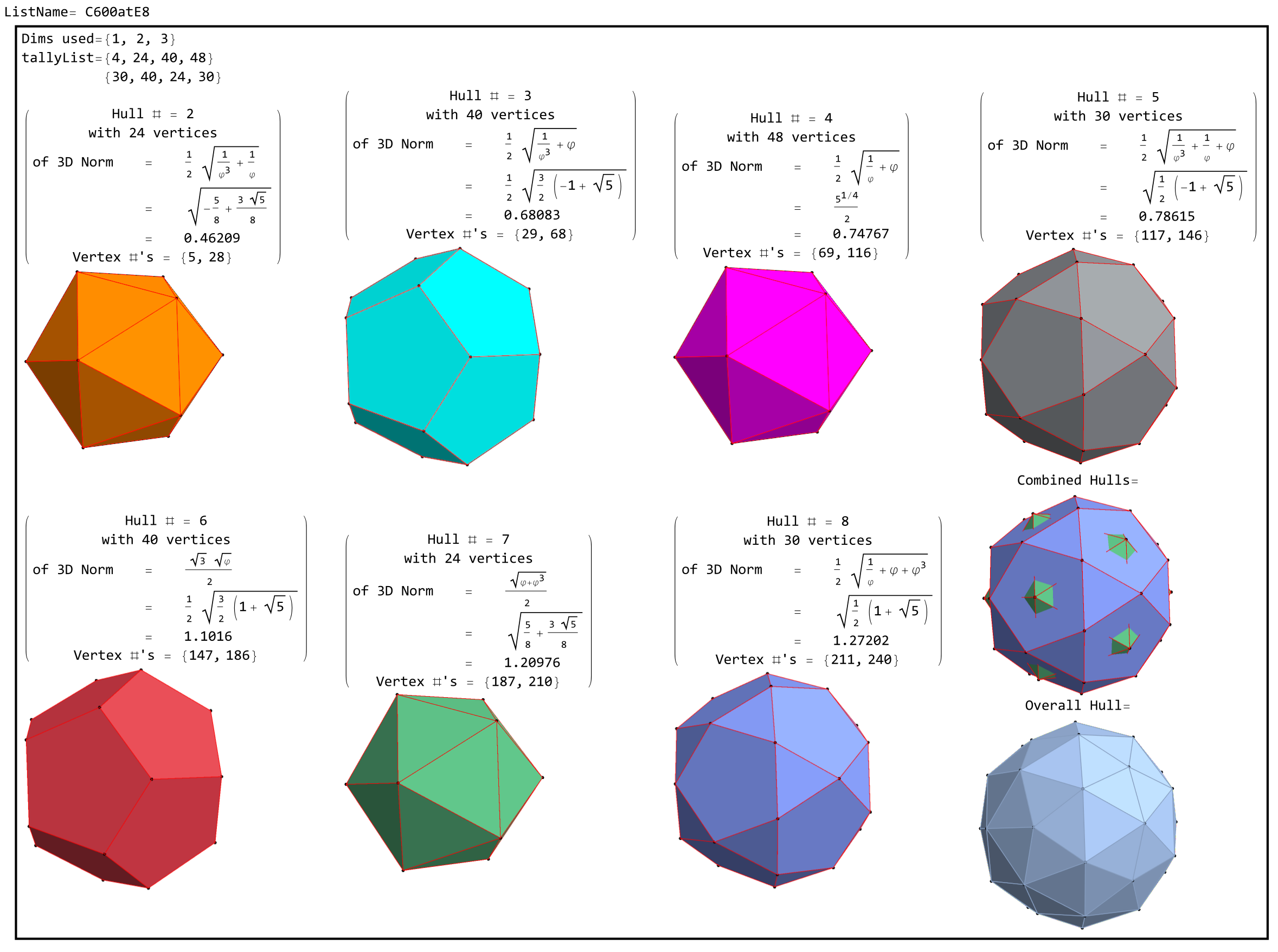}
\caption{\label{fig:421-concentric-hulls}Concentric hulls of $4_{21}$ in Platonic 3D projection with numeric and symbolic norm distances}
\end{figure}

\begin{figure}[!ht]
\center
\includegraphics[width=525pt]{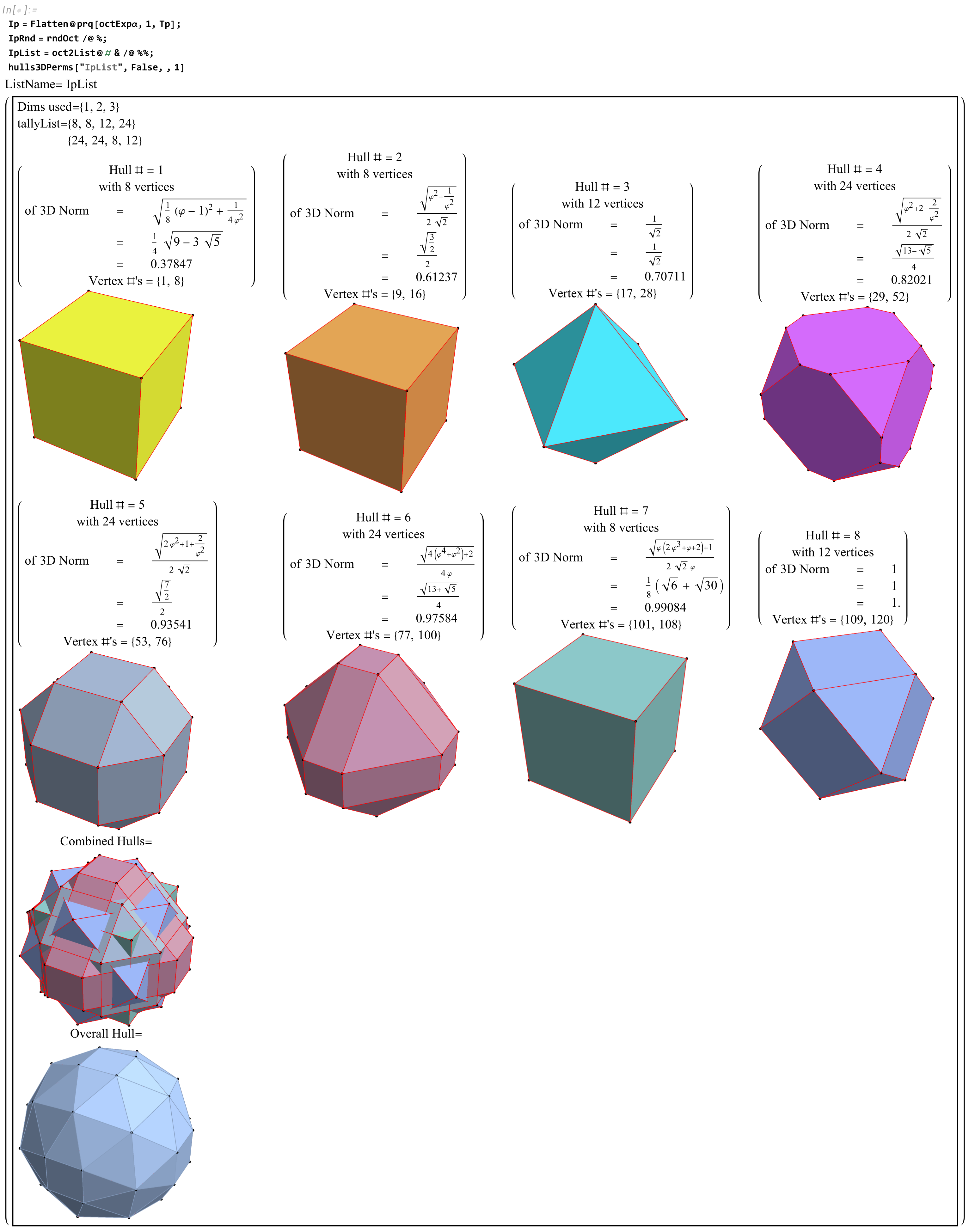}
\caption{\label{fig:Ip}Concentric hulls of I' as the parent $H_4$ 600-cell of order 120 in Platonic 3D projection with numeric and symbolic norm distances. This is generated by $\mathtt{I'=prq[}$$\alpha^{0-4}$$\mathtt{,1,T']}$.}
\end{figure}

\begin{figure}[!ht]
\center
\includegraphics[width=500pt]{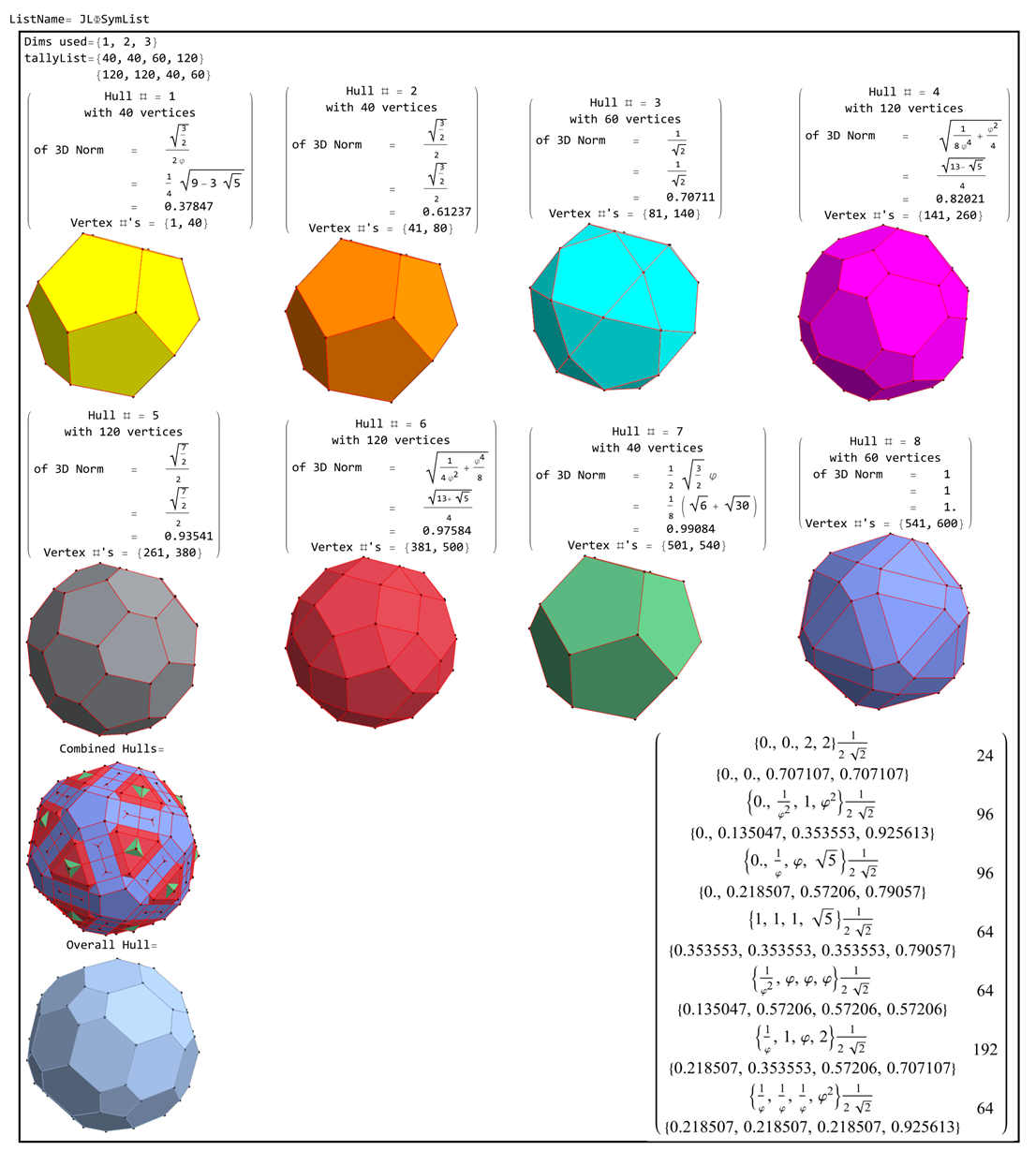}
\caption{\label{fig:J}Concentric hulls of J as the tri-rectified $H_4$ 120-cell of order 600 in Platonic 3D projection with numeric and symbolic norm distances. This is generated by $\mathtt{J=prq[A',1,I]=prq[A',}$$\alpha^{0-4}$$\mathtt{,T]}$.\\
Note: The numeric and symbolic tally list of unpermuted vertex values in the lower-right corner}
\end{figure}

\begin{figure}[!ht]
\center
\includegraphics[width=525pt]{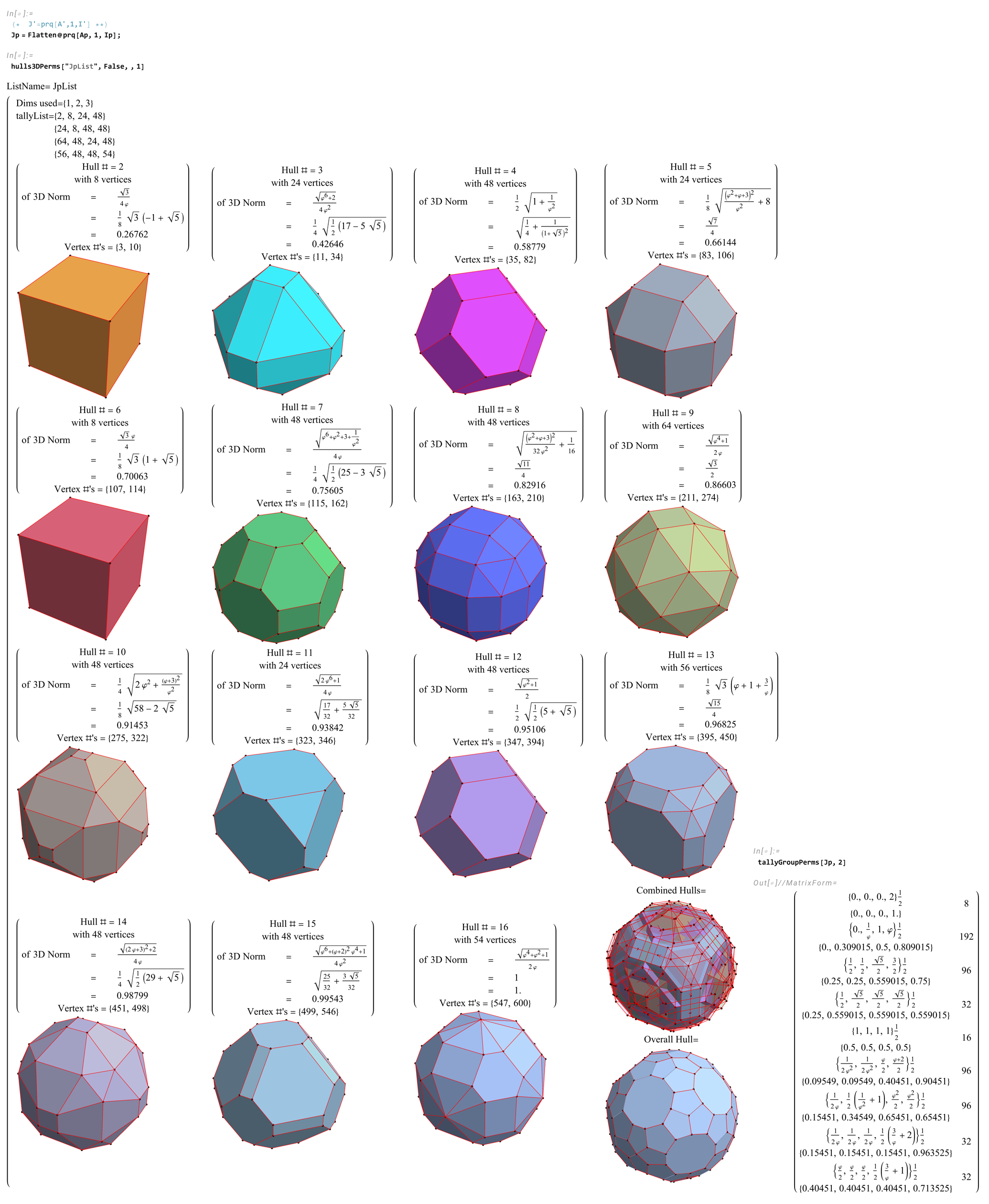}
\caption{\label{fig:Jp}Concentric hulls of J' as the tri-rectified $H_4$ 120-cell of order 600 in Platonic 3D projection with numeric and symbolic norm distances. This is generated by $\mathtt{J'=prq[A',1,I']=prq[A',}$$\beta^{0-4}$$\mathtt{,T']}$.\\
Note: The numeric and symbolic tally list of unpermuted vertex values in the lower-right corner}
\end{figure}

\section{\label{app:A}\textit{Concentric hulls from Platonic 3D projection with numeric and symbolic norm distances}\\
Figs. \ref{fig:421-concentric-hulls}-\ref{fig:Jp}
\ \\}

\begin{figure}[!ht]
\center
\includegraphics[width=525pt]{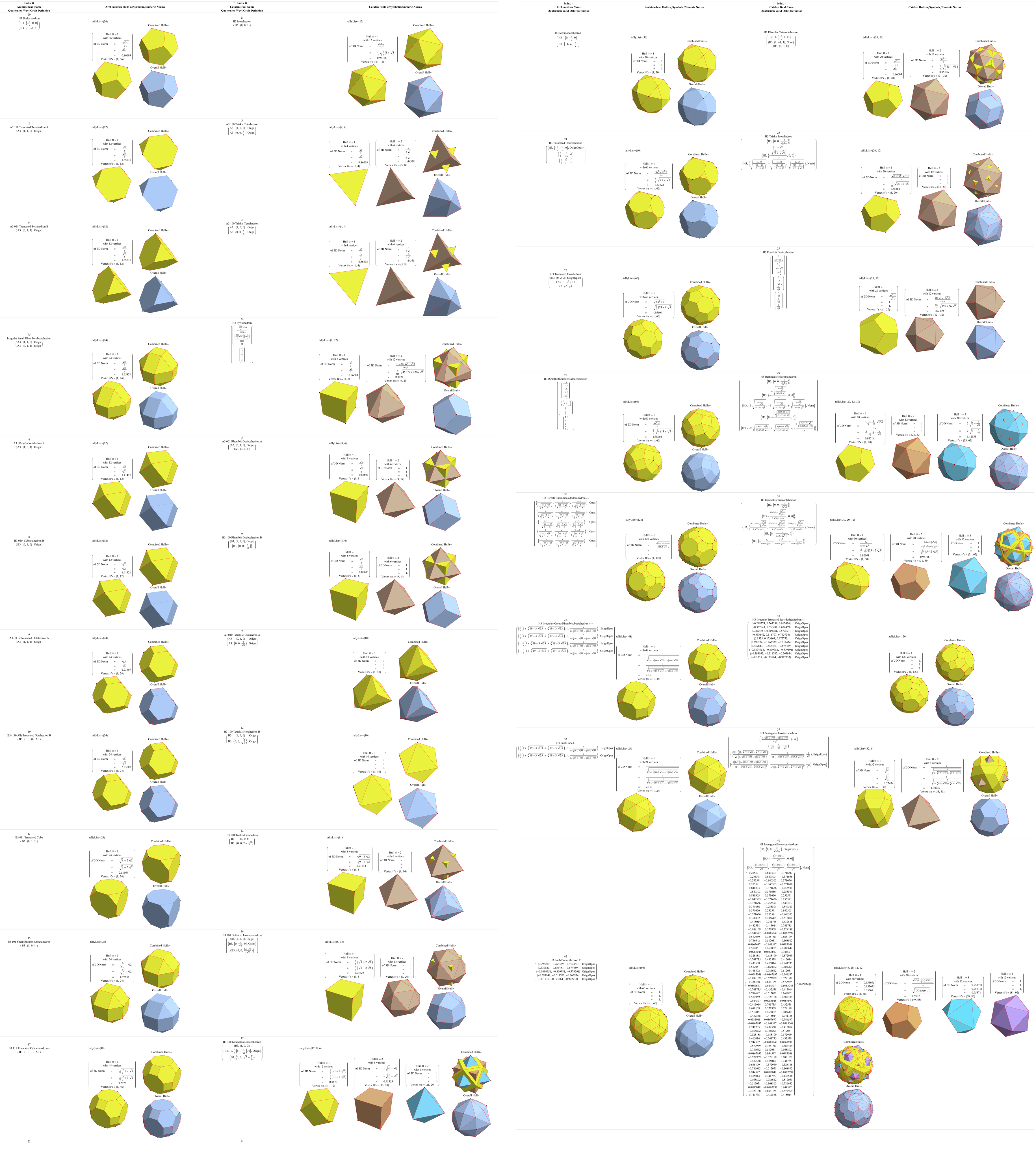}
\caption{\label{fig:Archimedean-and-Catalan-solids} \textit{Archimedean and dual Catalan solids, including their irregular and chiral forms. These were created using quaternion Weyl orbits directly from the $A_3$, $B_3$, and $H_3$ group symmetries\cite{koca2011chiral} listed in the first column.}}
\end{figure}

\section{\label{app:B}\textit{Archimedean and dual Catalan solids}\\
Fig. \ref{fig:Archimedean-and-Catalan-solids}
\ \\}

\begin{figure}[!ht]
\center
\includegraphics[width=525pt]{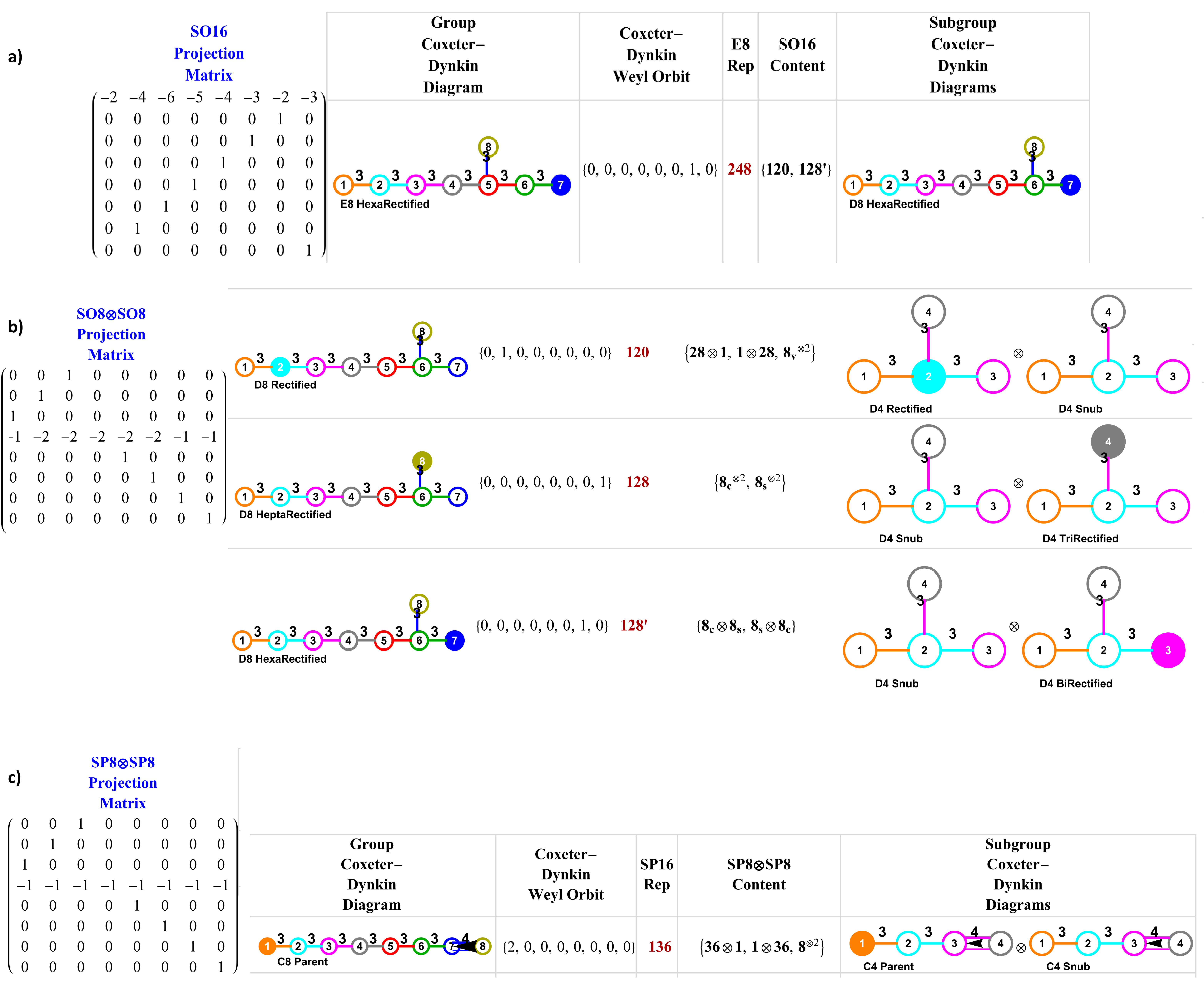}
\caption{\label{fig:CombinedEmbeddings}
Breakdown of $E_8$ maximal embeddings at height 248 of content SO(16)=$D_8$ (120,128')\\
a) Height 248 SO(16) content 120=(112+4+4)+128’\\
b) Height 120 and 128' SO(8)$\otimes$SO(8) content w/$8_{v,c,s}^{\otimes 2}$ triality\\
c) Height 136 Sp(8)$\otimes$Sp(8) content (32+4)$\otimes$1, 1$\otimes$(32+4), $8^{\otimes 2}$\\
Note: This output was created in \textit{Mathematica}\textsuperscript{TM} with support from the GroupMath\cite{Fonseca_2021} and SuperLie\cite{Grozman_2004} packages.}
\end{figure}

\begin{figure}[!ht]
\center
\includegraphics[width=525pt]{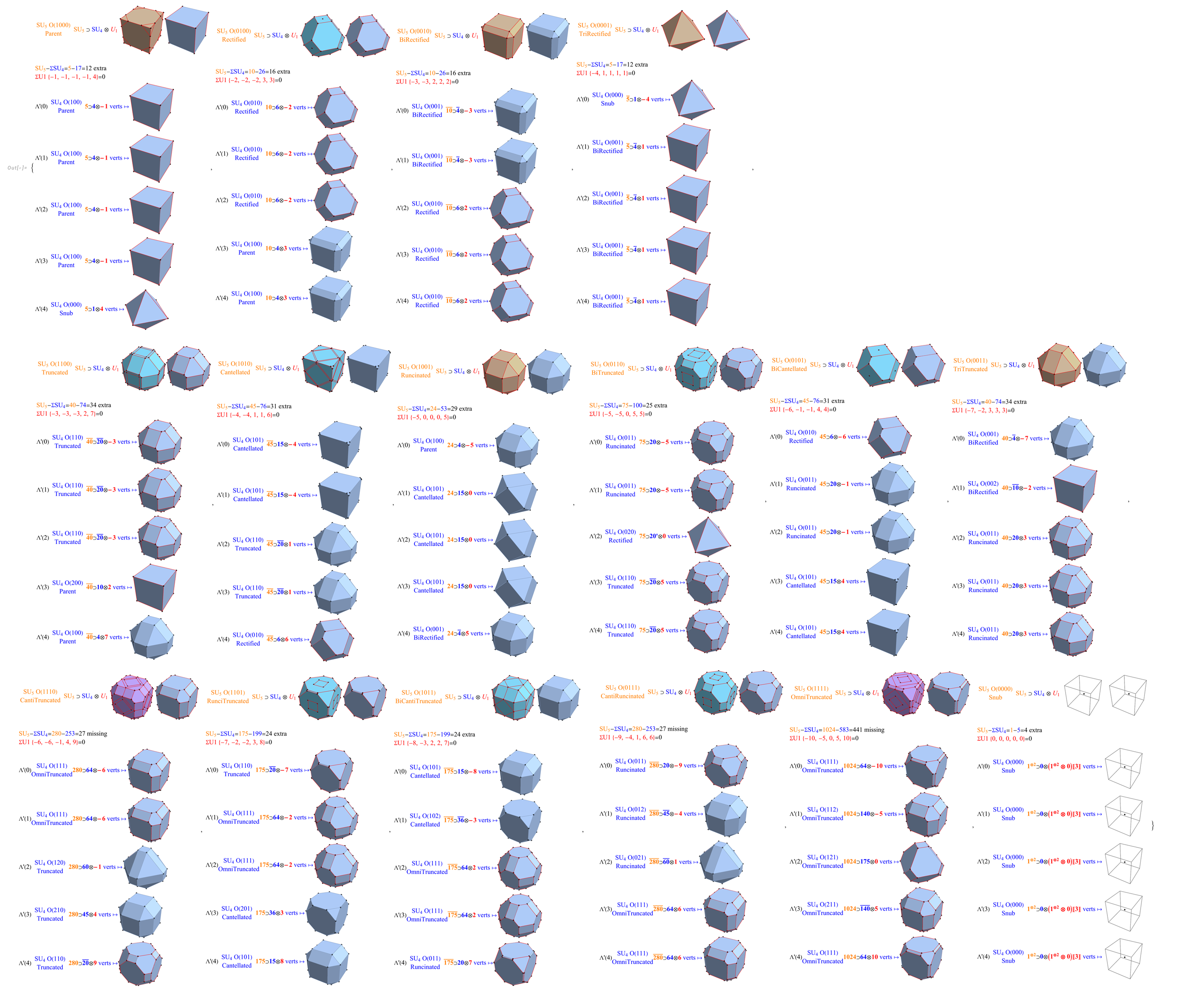}
\caption{\label{fig:A4-A3-SU5xSU4xU1-1024-Koca}$A_3$ in $A_4$ embeddings of SU(5)$\supset$SU(4)$\otimes$$U_1$\\
These include the specified 3D quaternion Weyl orbit hulls for each subgroup identified.}
\end{figure}

\section{\label{app:C}\textit{Maximal SO(16)=$D_8$ related embeddings of $E_8$ at height 248}\\
Figs. \ref{fig:CombinedEmbeddings}-\ref{fig:A4-A3-SU5xSU4xU1-1024-Koca}
\ \\}

\begin{figure}[!ht]
\center
\includegraphics[width=480pt]{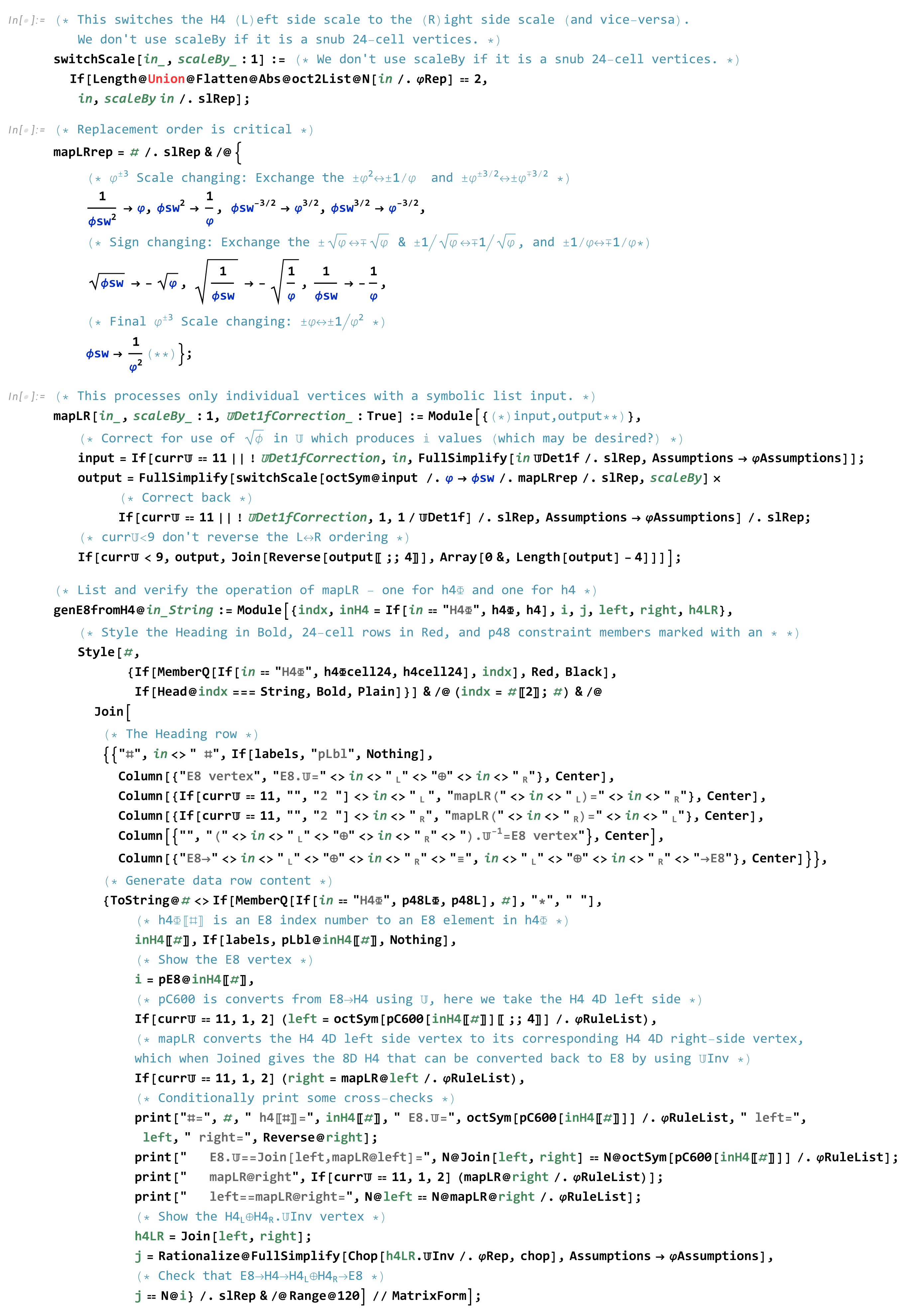}
\caption{\label{fig:Isomorphism-code} \textit{Mathematica}\textsuperscript{TM} code to generate the output showing $E_8\leftrightarrow H_4$ isomorphism}
\end{figure}

\begin{figure}[!ht]
\center
\includegraphics[width=525pt]{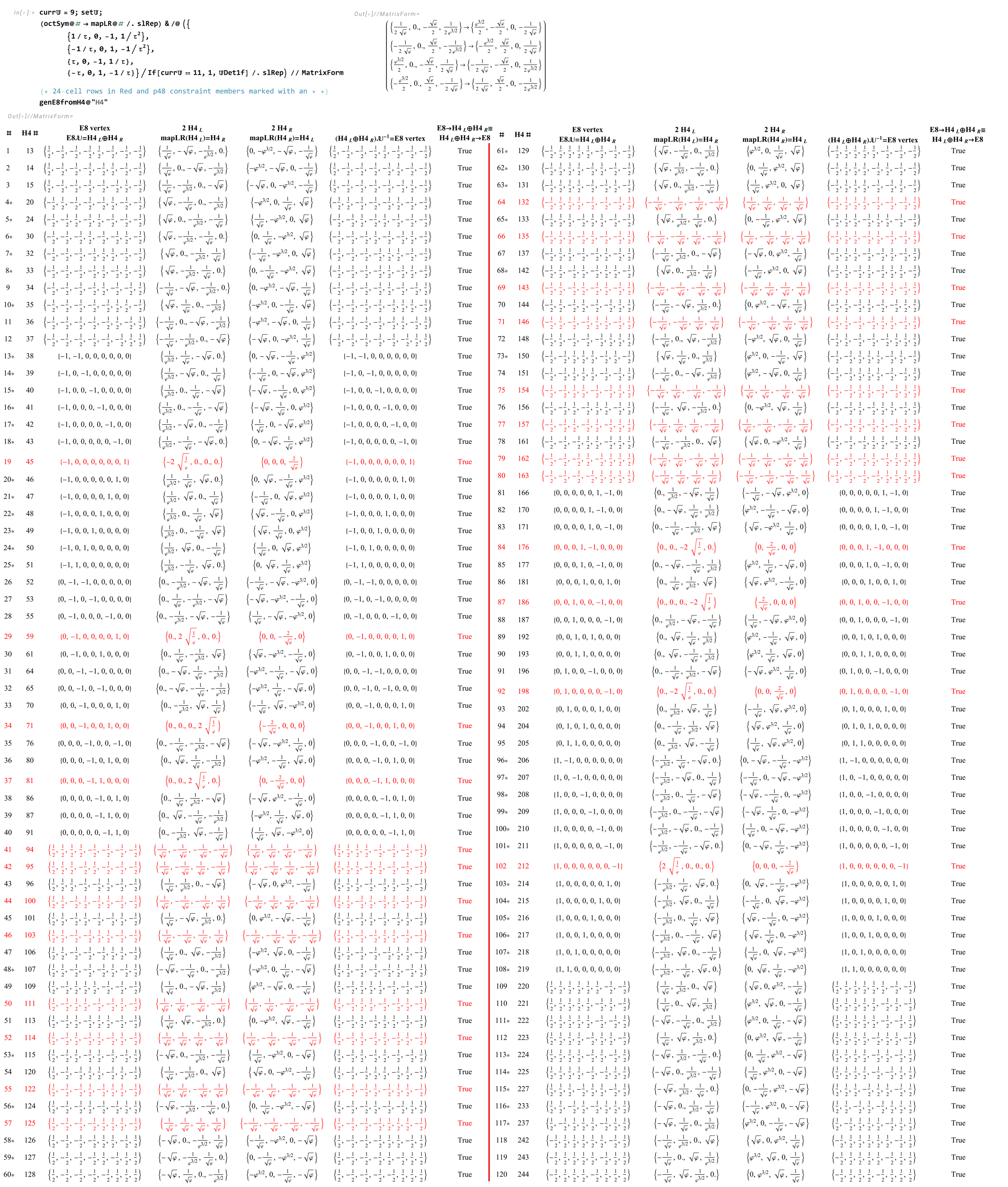}
\caption{\label{fig:Isomorphism-output-H4} Output showing detail of $E_8$$\leftrightarrow$$H_{4}$(L$\oplus$R) isomorphism for each vertex\\
Note: Red rows indicate $D_4$ 24-cell membership and the * identifies those satisfying the constraint of a unit normed p$\in$$S_L$ where $ p^0=\vert p^5\vert =\vert \bar{p}^5\vert =1\land \bar{p}^1=\pm p^4\land \bar{p}^4=\pm p\land \bar{p}^2=p^3\land \bar{p}^3=p^2$.}
\end{figure}

\begin{figure}[!ht]
\center
\includegraphics[width=525pt]{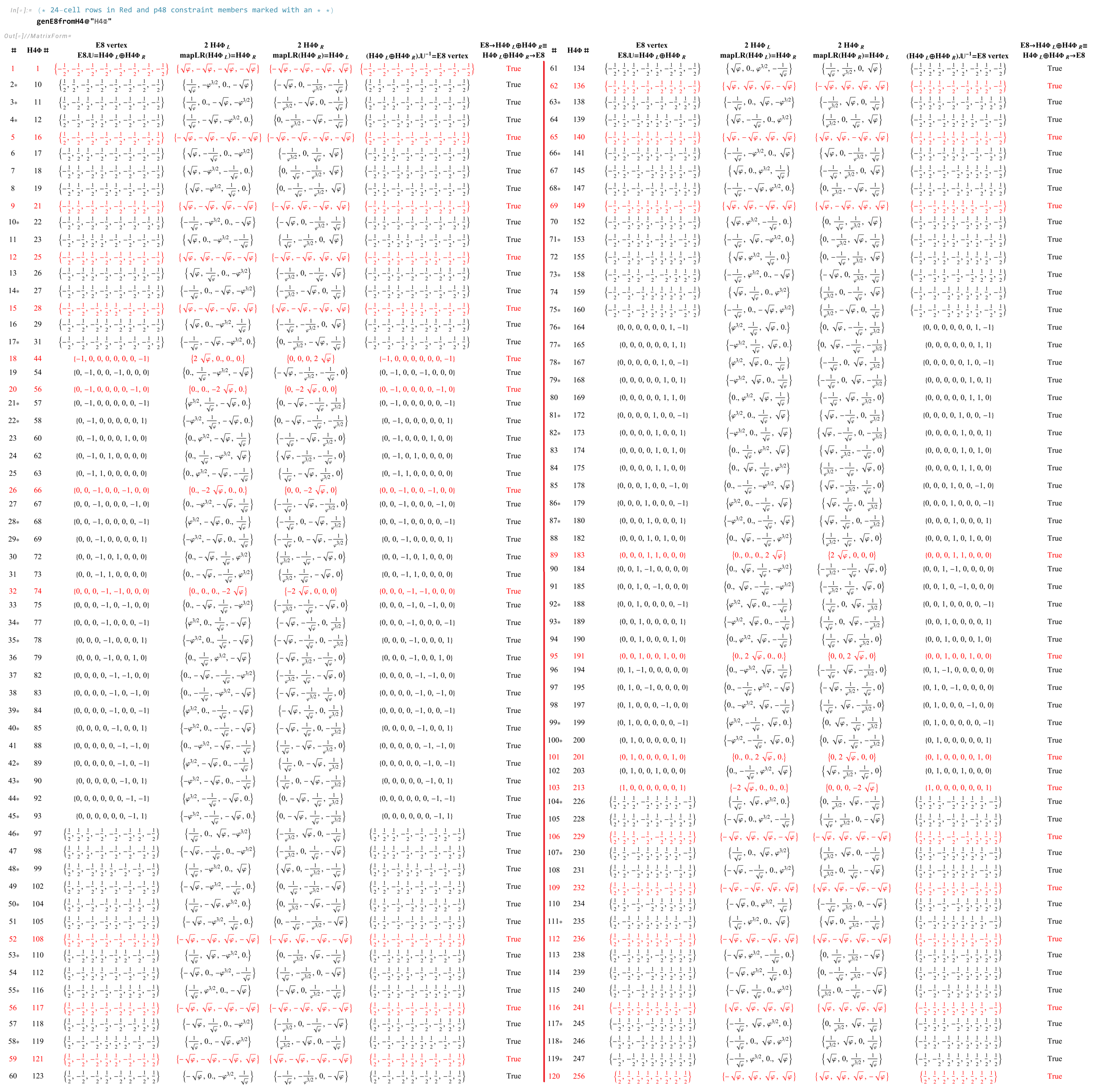}
\caption{\label{fig:Isomorphism-output-H4phi} Output showing detail of $E_8$$\leftrightarrow$$\varphi H_{4}$(L$\oplus$R) isomorphism for each vertex\\
Note: Red rows indicate $D_4$ 24-cell membership and the * identifies those satisfying the constraint of a unit normed p$\in$$\varphi S_L$ where $ p^0=\vert p^5\vert =\vert \bar{p}^5\vert =1\land \bar{p}^1=\pm p^4\land \bar{p}^4=\pm p\land \bar{p}^2=p^3\land \bar{p}^3=p^2$.}
\end{figure}

\section{\label{app:D}\textit{\textit{Mathematica}\textsuperscript{TM} code and output showing $E_8\leftrightarrow H_4$ isomorphism}\\
Figs. \ref{fig:Isomorphism-code}-\ref{fig:Isomorphism-output-H4phi}\ \\}

\end{document}